\documentclass[11pt,reqno]{amsart}
\usepackage{color, amsmath,amssymb, amsfonts, amstext,amsthm, latexsym}
\allowdisplaybreaks
\newcommand{\R}{{\mathbb R}}

\newcommand{\cF}{{\cal F}}
\newcommand{\RR}{{\mathbb R}}

  \newcommand{\EX}{{\mathbb E}}
  \newcommand{\PX}{{\mathbb P}}
\newcommand{\PP}{{\mathbb P}}

\newcommand{\nd}{{\nabla \cdot}}
\newcommand{\dn}{{\cdot \nabla}}

\newcommand{\s}{\sigma}
\newcommand{\e}{\varepsilon}
\renewcommand{\a}{\alpha}

\newcommand{\om}{\omega}
\newcommand{\Om}{\Omega}
\newcommand{\D}{\Delta}
\newcommand{\de}{\delta}

\newcommand{\pp}{\frac{\partial}{\partial t}}
\renewcommand{\k}{\kappa}
\renewcommand{\th}{\theta}

\renewcommand{\cF}{\mathcal F}

\numberwithin{equation}{section}
\newtheorem{theorem}{Theorem}[section]
\newtheorem{defn}[theorem]{Definition}

\newtheorem{lemma}[theorem]{Lemma}
\newtheorem{remark}[theorem]{Remark}
\newtheorem{prop}[theorem]{Proposition}

\begin{document}

\title[Large deviations for the Boussinesq model]
{Large deviations for the Boussinesq Equations  under Random
Influences}

\author[J. Duan and A. Millet]
{Jinqiao Duan and Annie Millet  }

\address[J.~Duan]
{Department of Applied Mathematics\\
Illinois Institute of Technology\\
Chicago, IL 60616, USA} \email[J.~Duan]{duan@iit.edu}

\address[A.~Millet]
{ SAMOS-MATISSE, Centre d'\'Economie de la Sorbonne (UMR 8174), Universit\'{e} Paris 1,
 Centre Pierre Mend\`{e}s France,
90 rue de Tolbiac, F- 75634 Paris Cedex 13, France {\it and}
Laboratoire de Probabilit\'es et Mod\`eles Al\'eatoires (UMR 7599),
       Universit\'es Paris~6-Paris~7, Bo\^{\i}te Courrier 188,
          4 place Jussieu, 75252 Paris Cedex 05, France } \email[A.
~Millet]{amillet@univ-paris1.fr {\it and} annie.millet@upmc.fr}

\thanks{This research was partly supported by the NSF Grant
   0620539 (J. Duan) and       by the research project BMF2003-01345 (A. Millet). }

\subjclass[2000]{Primary 60H15, 60F10; Secondary  35R60, 76D05, 76R05 }

\keywords{Boussinesq equations, B\'{e}nard convection,  large
deviations, stochastic PDEs, stochastic Navier-Stokes equations,
impact of noise on system evolution, multiplicative noise}

\begin{abstract}
A Boussinesq model for the B\'{e}nard convection under random
influences is considered as a system of stochastic partial
differential equations. This is a coupled system of  stochastic
Navier-Stokes equations and the transport equation for
temperature. Large deviations are proved, using a weak convergence
approach based on a variational representation of functionals of
infinite-dimensional Brownian motion.

\end{abstract}

\maketitle

\section{Introduction}\label{s1}

The need to take stochastic effects into account for modeling
complex systems has now become widely recognized. Stochastic
partial differential equations (SPDEs) arise naturally as mathematical
models for nonlinear macroscopic dynamics under random influences.
It is thus desirable to understand the impact of such random
influences on the system evolution \cite{WaymireDuan, PZ92, Roz}.

 The Navier-Stokes equations are often
coupled with other equations, especially, with the scalar
transport equations for   fluid density, salinity, or temperature.
These coupled equations (often with the Boussinesq approximation)
model a variety of phenomena in environmental, geophysical, and
climate systems \cite{DijkBook, DuanGaoSchm, Oz3}. We consider the
Boussinesq equations in which the scalar quantity is temperature,
under different   boundary conditions for the temperature at
different parts (top and bottom) of the boundary. This is a
B\'{e}nard convection problem. With other boundary conditions, the
Boussinesq equations model various phenomena in weather and
climate dynamics, for example.
We take random forcings into account and formulate the B\'{e}nard
convection problem as a system of stochastic partial differential
equations (SPDEs). This is a coupled system of the stochastic
Navier-Stokes equations and the stochastic transport equation for
temperature.

In various papers about large deviation principle (LDP) for
solutions $u^\e$
 to SPDEs or to evolution equations
 in a semi-linear framework
 \cite{Cerrai-Rockner, Chang, ChenalMillet,  Chow, FW92,
Kallian-Xiong, Peszat, Sowers, Zab},
the strategy used is similar to the classical one for diffusion processes.
 A very general version of Schilder's theorem yields the LDP for
the Gaussian noise $\sqrt{\e} W$ driving the stochastic forcing term,
 with a good rate function $\tilde{I}$ written in terms of its reproducing kernel
Hilbert space (RKHS).
 However, since the noise is not additive, the process  $u^\e$  is not a continuous function
of the noise, which creates technical difficulties. As if the contraction
principle were true, one defines deterministic controlled equations $u_h$ which are
similar to the stochastic one, replacing the stochastic integral with respect to
the noise $\sqrt{\e} W$ by deterministic integrals in terms of elements $h$ of its
RKHS. Once well-posedness of this controlled equation is achieved,
one proves that  solution $u^\e$ to the
 stochastic  evolution equation  satisfies a LDP with a rate function
  ${I}$  defined in terms of $\tilde{I}$
 and of $u_h$ by means of an energy minimization problem.
  In order to transfer the LDP from the noise to the process,
there are two classical proofs, each of which contains two main steps.
 One way consists in proving
 a continuity property of the map $h\mapsto u_h$ on level sets of the
rate function $\tilde{I}$ and then some Freidlin-Wentzell inequality,
which  states  continuity of the process with respect to the noise
except on an exponentially small set. Another classical method in
proving LDP for evolution equations is to establish both some exponential tightness and
exponentially  good approximations for some  approximating sequence where the diffusion
coefficient is stepwise constant. These methods require some time H\"older regularity  that one can
obtain when the diffusion coefficient is controlled in term of the $L^2$-norm
of the solution, but not in the framework we will use here, where
the bilinear term creates technical problems.
An alternative approach \cite{Feng} for large deviations is based
 on nonlinear semi-group theory and infinite-dimensional
 Hamilton-Jacobi equations, and it also requires  establishing 
 exponential tightness.

The method used in the present paper is related to the Laplace
principle. One proves directly that the level sets of the rate
function ${I}$ are compact and then establishes  weak convergence
of solutions to stochastic controlled equations written in terms
of the noise $\sqrt{\e} W$ shifted by a random element $h_\e$ of
its RKHS. This is again some kind of continuity property written
in terms of the distributions. Unlike \cite{Sundar},
well-posedness and a priori estimates are proved directly for very
general stochastic controlled equations with a forcing term including a stochastic
integral and a deterministic integral with respect to a random
element $h_\e$ of the RKHS of the noise, and for diffusion coefficients which
 may depend on the gradient.
 Indeed, if the
well-posedness for the stochastic controlled equation can be
deduced from that of the stochastic equation by means of a
Girsanov transformation, the a priori estimates uniform in $\e>0$,  which are a key
ingredient of the proof of the weak convergence result, cannot be
deduced from the corresponding ones for the stochastic B\'enard
equation since as $\e\to 0$, the $p>1$ moments of the Girsanov density go to
infinity exponentially fast.
 Well-posedness has been proved in \cite{Ferrario} for the
stochastic Boussinesq equation only in the particular case of an
additive noise on the velocity component. This weak convergence
approach has been introduced in \cite{BD00, BD07}. This  method
has been   recently applied to
  SPDEs \cite{Sundar, WangDuan} or SDEs   in infinite dimensions \cite{Ren}.
  Finally note that the proofs of the weak convergence and compactness
property require more assumptions
on the diffusion coefficient $\sigma$ which may not depend on the gradient.
 Indeed, in order to prove convergence of integrals
defined in terms of elements $h_\e$ of the RKHS  of the noise only using weak convergence of
$h_\e$, we also  need to deal with localized integral estimates of  time increments.
With additional assumptions on the diffusion coefficient
 we are  able to provide complete details of the proof  of this statement which
was missing in \cite{Sundar}.

This paper is organized as follows. The mathematical formulation
for the stochastic  B\'{e}nard model is in \S \ref{s2}. Then
the well-posedness  and  general a priori estimates for the model are proved  in \S \ref{s3}.
Finally, a large deviation principle is shown in \S \ref{s4}.

\section{Mathematical formulation} \label{s2}
Let $D =(0, l) \times (0, 1)$ be a rectangular   domain  in the
vertical plane. Denote by $x=(x_1,x_2)$ the spatial variable, $ u
= (u_1, u_2)$ the velocity field, $p$ the pressure field, $\th$
the temperature field, and $(e_1, e_2)$ the standard basis in
 $\R^2$. 

 We consider the following stochastic coupled Navier-Stokes and heat transport
 equations for the B\'{e}nard convection
 problem \cite{Foias}:
\begin{eqnarray}
\pp u^\e + u^\e \dn u^\e-\nu \D u^\e + \nabla p^\e &=&  \th^\e e_2  + \sqrt{\e}\; n_1(t),
\qquad \nd u^\e   = 0, \label{eqn1}\\
\pp \th^\e +u^\e \dn \th^\e -u_2^\e -\k\D \th^\e &=& \sqrt{\e} \;
n_2(t),\label{eqn3}
\end{eqnarray}
with boundary conditions
\begin{eqnarray}
u^\e =0\;\;  \& \;\; \th^\e=0 \;\; \mbox{on}\;\;x_2=0\; \mbox{and}\;x_2=1,   \\
u^\e, p^\e, \th^\e, u_{x_1}^\e, \th_{x_1}^\e \; \mbox{are periodic in}\; x_1 \;
\mbox{with period}\; l,
\end{eqnarray}
where $n_1, n_2$ are noise forcing terms and $ \e>0 $ is a small
parameter.
\medskip

We consider the abstract functional setting for this system as in
\cite{Foias, Ferrario}; see also \cite{Constantin, Temam}.
Let $L^2(D)$ be endowed with the usual scalar product and the
induced norm. Consider another Hilbert space of vector-valued functions:
\begin{align*}
  \dot{\textbf{L}}^2(D) =& \{u\in L^2(D)^2,\; \nd u=0, \; u |_{x_2=0}=u |_{x_2=1}=0,\;
u \; \mbox{is periodic in}\; x_1 \; \mbox{with period}\; l  \}    \\
  \dot{L}^2(D) =& \{\th \in L^2(D),\;  \; \th|_{x_2=0}=\th|_{x_2=1}=0,\;
\th \; \mbox{is periodic in}\; x_1 \; \mbox{with period}\; l  \}
\end{align*}
Let $H= \dot{\textbf{L}}^2(D) \times \dot{L}^2(D)$
be the product Hilbert space.
 We   denote by the same notations, $(\cdot, \cdot)$ and $|\cdot |$,
  the scalar
 product and the induced norm, in  $\dot{\textbf{L}}^2(D)$, $\dot{L}^2(D)$
 and  $H$,
  $$ (\phi, \psi)=\int_D \phi(x) \psi(x)dx,
   \; \; |\phi|=\sqrt{(\phi, \phi)}=\sqrt{|\phi_1|^2 + |\phi_2|^2}.
  $$
 Define $V=V_1 \times V_2$, where
\begin{align*}
 V_1=&\{v\in H^1(D)^2: \nd v=0,\; v|_{x_2=0}=v|_{x_2=1}=0; \; v
\; \mbox{is periodic in}\; x_1 \; \mbox{with period}\; l \}, \\  
  V_2=&\{f\in
H^1(D): f|_{x_2=0}=f|_{x_2=1}=0; \; f \; \mbox{is periodic in}\;
x_1 \; \mbox{with period}\; l \}.
\end{align*}
 Then  $V$  is a product
Hilbert space with the scalar product and the induced norm,
$$
((\phi, \psi))=\int_D \nabla \phi \cdot \nabla \psi dx, \;\;
\|\phi\|=\sqrt{((\phi, \phi))} =\sqrt{\|\phi_1\|^2 +
\|\phi_2\|^2},
$$
where, 
to ease the notation,  the space variable $x$ is omitted when
writing integrals on $D$.  Again, we also use  the same notations
for
  the scalar
 product and the induced norm  in  $V_1$ and $V_2$.
Let $V'$ be the dual space of $V$. We have the dense and
continuous embeddings
$ V \hookrightarrow H=H' \hookrightarrow V' $
and denote by $\langle \phi , \psi\rangle$ the duality between $\phi\in V$ (resp. $V_i$) and
$\phi\in V'$ (resp. $V'_i$).
Recall that there exists some positive constant $c_1$ such that for $u\in V_1$, $\theta\in V_2$,
\begin{equation}\label{VL4}
|u|_{L^4(D)^2}^2 \leq c_1\, |u|\, \|u\|, \quad \mbox{\rm and}\quad |\theta|_{L^4(D)}^2 \leq c_1\, |\theta|\,
\|\theta\|.
\end{equation}
Furthermore, the Poincar\'e inequality yields the existence of a positive constant $c_2$
such that
\begin{equation} \label{Poincare}
|\phi|\leq c_2\, \|\phi\|, \quad \forall \phi\in V.
\end{equation}
To lighten the  notations, we will set for $\phi=(u,\theta)$,  $u\in L^4$, $\theta\in L^4$ and
$\phi\in L^4$ for vectors of dimensions  2,1 and 3 whose components belong to $L^4(D)$
 and denote the corresponding norms by $|\, |_{L^4}$.

Consider an unbounded linear operator $A=( \nu A_1, \k A_2 ): H \to
H$ with $D(A)=D(A_1) \times D(A_2)$ where $D(A_1)=V_1\cap
H^2(D)^2$, $D(A_2)=V_2 \cap H^2(D)$ and define
\begin{eqnarray*}
\langle A_1 u,v\rangle =((u, v)), \; \langle A_2 \th,\eta\rangle =((\th, \eta)),
  \; \forall u, v \in D(A_1), 
\; \forall \th, \eta \in D(A_2).
\end{eqnarray*}
Both the Stokes operator $A_1$ and the Laplace operator $A_2$ are
self-adjoint, positive, with compact self-adjoint inverses.
 They map $V$ to $V'$.
We also introduce the bilinear operators $B_1$ and $B_2$ as
follows: for $ u, v, w \in V_1$ and $ \th, \eta \in V_2$,
\begin{eqnarray*}
\langle B_1(u,v),w\rangle  &=& \int_D [u\dn v]wdx \quad
 := \sum_{i,j=1,2} \int_D u_i \, \partial_i v_j\,  w_j dx,  \\
\langle B_2(u,\th),\eta\rangle &=&\int_D[u\dn\th] \eta dx\quad
 := \sum_{i=1,2} \int_D u_i \, \partial_i \,  \theta\;  \eta \, dx .
\end{eqnarray*}
With the notation $\phi^\e=( u^\e,\th^\e)$
 and under the above
formulation, we assume that the noise terms $n_1$ and $n_2$ are
 respectively
$ \s_1(t,\phi) \pp W^1(t)$ , $ \s_2(t,\phi) \pp W^2(t)$,
 where  $W^1(t), W^2(t)$ are independent  Wiener processes  defined   on a filtered
probability space $(\Om, \cF, \cF_t, \PX)$, taking values in
$\dot{\textbf{L}}^2(D)$ and
   $\dot{L}^2(D)$,   with linear symmetric
positive covariant operators $Q_1$ and $Q_2$,  respectively. We
denote $Q =(Q_1, Q_2) $.
 It is a linear symmetric positive
covariant operator in the Hilbert space $H$. We assume that $Q_1,
Q_2$ and thus $Q$ are trace class (and hence compact \cite{PZ92}),
i.e., $tr(Q) < \infty$.

As in \cite{Sundar}, let $H_0 = Q^{\frac12} H$. Then $H_0$ is a
Hilbert space with the scalar product
$$
(\phi, \psi)_0 = (Q^{-\frac12}\phi, Q^{-\frac12}\psi),\; \forall
\phi, \psi \in H_0
$$
together with the induced norm $|\cdot|_0=\sqrt{(\cdot,
\cdot)_0}$. The embedding $i: H_0 \to  H$ is Hilbert-Schmidt and
hence compact, and moreover, $i \; i^* =Q$.

Let $L_Q$ be the space of linear operators $S$ such that
$SQ^{\frac12}$ is a Hilbert-Schmidt operator (and thus a compact
operator \cite{PZ92}) from $H$ to $H$. The norm in the space $L_Q$ is
  defined by
 $|S|_{L_Q}^2 =tr (SQS^*)$,
 where $S^*$ is the adjoint operator of
$S$.

Note that the above formulation is equivalent to projecting the
first governing equation from   $\dot{L}^2(D)^2$   into the
``divergence-free" space      
 and thus the pressure term is absent.
 With these notation, the above Boussinesq system
(\ref{eqn1})-(\ref{eqn3}) becomes
\begin{eqnarray}
d u^\e + [\nu A_1 u^\e +B_1(u^\e,u^\e)-  \th^\e e_2]dt&=& \sqrt{\e} \;  \s_1(t,\phi^\e) \; dW^1(t),
 \label{benard1}\\
d \th^\e +[\k A_2\th^\e +B_2(u^\e,\th^\e)-u_2^\e]dt &=& \sqrt{\e} \;  \s_2(t,\phi^\e)\;
dW^2(t).\label{benard2}
\end{eqnarray}
Thus,  we write this system  for $\phi^\e=(u^\e,\theta^\e)$ as
\begin{eqnarray}\label{Benard}
d \phi^\e + [A\phi^\e +B(\phi^\e)+R\phi^\e]dt = \sqrt{\e} \; \s(t, \phi^\e) dW(t),
\quad \phi^\e(0)=\xi:= ( u_0^\e, \th_0^\e),
\end{eqnarray}
where  $ W(t) = (W^1(t), W^2(t))$ and
\begin{eqnarray}
A\phi &=&( \nu A_1 u,\k A_2\th ), \\
 B(\phi) &=& (B_1(u, u), B_2(u, \th) ),  \\
  R\phi &=&(- \th e_2, -u_2), \\
  \s(t,\phi) &=&  (\s_1(t,\phi), \s_2(t,\phi))  . 
\end{eqnarray}
The noise intensity $\s: [0, T]\times V \to L_Q(H_0, H)$ is
assumed to satisfy the following: 
 \medskip

\noindent \textbf{Assumption A:}
There exist positive constants $K$ and $L$
such that\\
{\bf (A.1)}  $\s \in C\big([0, T] \times H; L_Q(H_0, H)\big)$  \\
{\bf (A.2)}   $|\s(t,\phi)|^2_{L_Q} \leq K (1+\|\phi\|^2), \quad \forall t\in [0,T]$ ,
$\forall \phi\in V$.\\
{\bf (A.3)}   $|\s(t,\phi)-\s(t,\psi)|^2_{L_Q} \leq L \|\phi -\psi\|^2,\quad \forall t\in [0,T]$,
$\forall \phi, \psi \in V$.
\medskip

  In what follows, to ease the notation,  we will suppose that $\sigma(t,\phi)=\sigma(\phi)$;
 however, all the
results have a straightforward extension to
time-dependent noise intensity under the assumption A.
When no confusion arises, we set $L^p:=L^p(D)$ for
$1\leq p<+\infty$ and denote by $C$ a constant which may change from one line to the next one.

\section{Well-posedness}  \label{s3}

The goal for this paper is to show the large deviation
principle 
 for  $(\phi^\e , \e>0)$ as $\e\to 0$, where $\phi^\e$ denotes the solution to the stochastic B\'enard equation 
(\ref{Benard}).

Let $\mathcal{A}$ be the  class of $H_0-$valued
$(\cF_t)-$predictable stochastic processes $\phi$ with the property
$\int_0^T |\phi(s)|^2_0 ds < \infty, \; $ a.s.
Let
\[S_M=\Big\{h \in L^2(0, T; H_0): \int_0^T |h(s)|^2_0 ds \leq M\Big\}.\]
The set $S_M$ endowed with the following weak topology is a
  Polish space (complete separable metric space)
\cite{BD07}:
$ d_1(h, k)=\sum_{i=1}^{\infty} \frac1{2^i} \big|\int_0^T \big(h(s)-k(s),
\tilde{e}_i(s)\big)_0 ds \big|,$
where $
\{\tilde{e}_i(s)\}_{i=1}^{\infty}$ is a complete orthonormal basis
for $L^2(0, T; H_0)$.
Define
\begin{equation} \label{AM}
 \mathcal{A}_M=\{\phi\in \mathcal{A}: \phi(\om) \in
 S_M, \; a.s.\}.
\end{equation}
 As in \cite{Sundar},  we prove existence and uniqueness of the solution to the B\'enard
equation. However,  in what follows, we will need some precise
bounds on the norm of the solution to a more general equation,
which contains  an extra forcing (or control) term driven by an
element of
    $ \mathcal{A}_M$. These required estimates cannot be deduced from the corresponding
ones by means of a Girsanov transformation. More precisely, let
  $h\in \mathcal{A}$, $\e\geq 0$
 and consider the following
generalized B\'enard equation with initial condition $\phi_h^\e(0)=\xi$.
For technical reasons, we need to add some control in the forcing term, with intensity
$\tilde{\s}\in C([0,T]\times H; L_Q(H_0,H))$ satisfying similar stronger conditions:\\
{\bf Assumptions $\tilde{\rm \bf A}$: } There exist positive constants $\tilde{K}$ and $\tilde{L}$ such that: \\
{\bf($\tilde{\rm\bf A}.1$)} \quad  $|\tilde{\s}(t,\phi)|_{L_Q}^2 \leq \tilde{K} (1+|\phi|_{L^4}^2)$,\quad
 $\forall t\in [0,T]$,
$\forall \phi\in L^4(D)^3$.\\
{\bf($\tilde{\rm\bf A}.2$)} \quad
 $|\tilde{\s}(t,\phi)-\tilde{\s}(t,\psi)|_{L_Q}^2 \leq \tilde{L} |\phi - \psi|_{L^4}^2$, \quad
 $\forall t\in [0,T]$,  $\forall \phi, \psi \in L^4(D)^3$.\\
Notice that since $V\subset L^4(D)^3$, the assumption  {\bf  $\tilde{\rm \bf A}$} is stronger
 than {\bf A}.
  For $\s$ ,
$\tilde{\s} \in C(H; L_Q(H_0,H))$ which  satisfy Assumptions  {\bf A} and  {\bf$\tilde{\rm\bf  A}$}
respectively,  set
\begin{equation} \label{Benardgene}
d \phi_h^\e(t) + \big[ A\phi_h^\e(t) +B(\phi_h^\e(t))+R\phi_h^\e(t)\big]dt =
\sqrt{\e} \s(\phi_h^\e(t)) dW(t) + \tilde{\s}(\phi_h^\e(t)) h(t) dt . 
\end{equation}
Recall that a stochastic process $\phi_h^\e(t,\om)$ is called the weak
solution for the generalized stochastic B\'{e}nard problem (\ref{Benardgene}) on $[0, T]$
 with initial condition $\xi$ if
$\phi_h^\e$  is in $C([0, T]; H) \cap L^2((0, T); V)$, a.s., and
satisfies
\begin{eqnarray}\label{weak}
 (\phi_h^\e(t), \psi)-(\xi, \psi) + \int_0^t \big[ (  \phi_h^\e(s), A\psi )
 + \big\langle  B(\phi_h^\e(s)), \psi\big\rangle + (R\phi_h^\e(s),\psi)\big]ds \nonumber \\
 =  \sqrt{\e} \int_0^t \big(\s(\phi_h^\e(s)) dW(s),\psi\big) +
  \int_0^t \big( \tilde{\s}(\phi_h^\e(s)) h(s)
\, ,\, \psi\big)\, ds,\;\; a.s.,
\end{eqnarray}
for all $\psi \in D(A)$ and all $t \in [0,T]$.
 In most of the
analysis here,  we work in the Banach space
$ X: = C\big([0, T]; H\big) \cap L^2\big((0, T);V\big) $
with the norm
\begin{eqnarray}\label{norm}
\|\phi\|_X = \Big\{\sup_{0\leq s\leq T}|\phi(s)|^2+ \int_0^T \|
\phi(s)\|^2 ds\Big\}^\frac12.
\end{eqnarray}
\begin{theorem}\label{wellposeness} (Well-Posedness and A priori bounds)\\
Fix $M>0$; then there exists $\e_0:=\e_0(\nu,\kappa, K, L, \tilde{K}, \tilde{L}, T, M)>0$,
 such that the following existence
and uniqueness result is true  for $0\leq \e\leq \e_0$.
Let the initial datum satisfy  $\EX |\xi|^4 <
\infty$, let  $h\in   \mathcal{A}_M$ and $\e \in [0,\e_0]$; then
there exists a pathwise unique weak solution $\phi_h^\e$ of the
generalized stochastic B\'{e}nard problem \eqref{Benardgene} with
initial condition  $\phi_h^\e(0)=\xi \in H$ and  such that $\phi_h^\e\in X$
a.s.
 Furthermore, there exists a constant $C_1:=C_1(  \nu,\kappa,K,L, T,M)$ such that
for $\e\in [0,\e_0]$ and $h\in {\mathcal A}_M$,
\begin{equation} \label{boundgeneral}
E\|\phi^\e_h\|_X^2 \leq 1+ E\Big( \sup_{0\leq t\leq T}
 |\phi_h^\e(t)|^4
+ \int_0^T \|\phi_h^\e(t)\|^2\, dt \Big) \leq C_1\, \big( 1+E|\xi|^4\big).
\end{equation}
 \end{theorem}
\begin{remark}
Note that if  $\sigma = 0$, i.e., when the noise term is absent, we
 deduce the existence and uniqueness of  the solution to the
``deterministic"
 control equation defined in terms of an element $h\in  L^2((0,T);H_0)$
 and an  initial condition $\xi \in H$
\begin{equation}\label{dc}
d \phi(t) + \big[ A\phi(t) +B(\phi(t))+R\phi(t)\big]dt =\tilde{\s}(\phi(t)) h(t) dt, \;\;
\phi(0)=\xi.
\end{equation}
 If  $h \in S_M$,  the solution $\phi$ to \eqref{dc} satisfies
\begin{eqnarray} \label{control-norm2}
 \sup_{0\leq s \leq T} |\phi(s) |^2+\int_0^T \|\phi(s)\|^2 ds
  \leq \tilde{C}_1(\nu, \k, \tilde{K}, \tilde{L}, T, M,  |\xi|).
\end{eqnarray}
\end{remark}
\begin{remark}
Finally, note that when $\phi_h^\e$ is a solution to the stochastic Boussinesq  equation
\eqref{Benard}, a similar argument shows that  Theorem \ref{wellposeness}
holds for any  $\e \geq 0$ if the  coefficients  $\sigma$ (resp.  $\tilde{\sigma}$) belong to
$C( [0,T]\times H ;  L_Q(H_0,H))$
and are such that in the upper estimates of
 the $L_Q$-norm appearing
in the right-hand sides of  conditions ({\bf A.2}) and ({\bf A.3})
 (resp. ($\tilde{\rm \bf A}.1$) and  ($\tilde{\rm\bf A}.2$),
one replaces the $V$ (resp. the $L^4$)  norms  of  $\phi$ and $\phi-\psi$ by their $H$-norms.

Indeed, in that case, for any fixed $\e>0$, the control of the $V$-norm of the solution,
or of its finite-dimensional approximation, only comes from the operators $A$ and $B$. Thus Lemmas
\ref{BB} and \ref{diffB} below prove that for $\alpha$ small enough, the $V$-norm can be dealt with.
\end{remark}
The proof of this theorem will require several steps.
The following lemmas gather some properties of $B_1$ and $B_2$.
 We refer  the reader to \cite{Constantin} or \cite{Temam} for
the  results on $B_1$ which are classical and sketch some  proofs of the corresponding results on $B_2$.
\begin{lemma} \label{lemBi}
For $u,v, w\in V_1$ and $\theta,\eta\in V_2$,
\begin{eqnarray*}
\langle B_1(u,v),v\rangle  &=&   0,    \quad
 \langle B_2(u,\th),\th\rangle \;  =\;  0,  \\ 
\langle B_1(u,v),w\rangle  &=&   - \langle B_1(u,w),v\rangle ,\quad  
\langle B_2(u,\th),\eta\rangle  \;=\;  -\langle B_2(u,\eta),\th\rangle . 
\end{eqnarray*}
\end{lemma}
Let $u\in  V_1$, $\theta\in V_2$ and $\phi=(u,\theta)\in V$; note
that $|\phi|^2=|u|^2+|\theta|^2$ and $\|\phi\|^2=\|u\|^2 +
\|\theta\|^2$. The following lemma provides
 upper bound estimates of
 $B_1$ and $B_2$.

\begin{lemma} \label{B2}
Let  $c_1$ denote the constant in \eqref{VL4}; then for any  $u\in V_1$,  
$\theta,\eta\in V_2$ and $\phi=(u,\theta)$,
 one has
\begin{eqnarray}
   |B_1(u,u)|_{V'_1}
  &\leq& |u|^2_{L^4} \leq  c_1\, |u| \;
\|u\|,  \label{normB1V'}    \\
|\langle B_2(u, \th), \eta\rangle | &\leq&
|u|_{L^4}\, |\th|_{L^4} \, \|\eta\|  \leq  c_1\, |\phi| \; \|\phi\| \|\eta\|.
 \label{inegB2}
\end{eqnarray}
\end{lemma}
\begin{proof}
We only check the properties on $B_2$. 
For $\phi=(u,\theta)\in V$ and $\eta\in V_2$,
Lemma \ref{lemBi}, 
 H\"older's inequality, and 
\eqref{VL4}  imply
\[ \big|\big\langle  B_2(u,\theta), \eta\big\rangle \big| =  \big|\big\langle  B_2(u,\eta), \theta\big\rangle\big|
 \leq \|\eta\|\, |u|_{L^4}\, |\theta|_{L^4} \leq c_1 \|\eta\|\, |u|^{\frac{1}{2}} \|u\|^{\frac{1}{2}}
|\theta|^{\frac{1}{2}} \|\theta\|^{\frac{1}{2}}. \] This yields
\eqref{inegB2}.
\end{proof}

\begin{lemma} \label{BB}
Let $\phi=(u,\theta) \in V$,  and let $v\in L^4(D)^2$ and $\eta\in L^4(D)$.
 For any   constant $\a >0$, the following estimates
hold:
 \begin{eqnarray}
|\langle B_1(u,u),v\rangle | &\leq &  \a \, \|u\|^2
+ \frac{3^3\, c_1^2 }{4^4\a^3}\; |u|^2 \; |v|^4_{L^4},  \label{BB1}    \\
 |\langle B_2(\phi), \eta\rangle| &\leq& \a\, \|\phi\|^2
+ \frac{3^3\, c_1^2}{4^4\a^3} \; |u|^2 \; |\eta|^4_{L^4}. \label{BB2}
\end{eqnarray}
\end{lemma}
\begin{proof}
We only  check \eqref{BB2}. The first part of
\eqref{inegB2}
and  Young's inequality yield 
\begin{align*}
 |\langle B_2(\phi), \eta\rangle|
 = |\langle B_2(u,\eta)\, ,\theta\rangle |
\, & \, \leq  |\eta|_{L^4} \;|u|_{L^4} \;
  |\nabla\th|_{L^2} 
\, \leq  \, \sqrt{c_1}\,    |\eta|_{L^4}\; |u|_{L^2}^\frac12 \;|\nabla
u|_{L^2}^\frac12 \;  |\nabla\th|_{L^2} \\
& \leq  \sqrt{c_1}\,  |\eta|_{L^4}\; |u|_{L^2}^\frac12 \;\|\phi\|^\frac32  
\; \leq \; \a\,  \|\phi\|^2 + \frac{3^3\, c_1^2}{4^4\a^3} \; |u|^2 \;
|\eta|^4_{L^4}.
\end{align*}
\end{proof}
The following lemma allows  rewriting differences of $B_i$ for $i=1,2$ and  deducing 
estimates for the difference of  $B$.
\begin{lemma} \label{diffB}
Let $\phi=(u,\theta)$ and $\psi=(v,\eta)$ belong to $V$. Then
\begin{eqnarray*}
 \big\langle  B_1(u,u)-B_1(v,v),u-v\big\rangle &= &- \big\langle  B_1(u-v,u-v),v\big\rangle ,\\
 \big\langle  B_2(\phi)-B_2(\psi), \theta-\eta \big\rangle
&= &- \big\langle B_2(\phi-\psi),\eta\big\rangle .
\end{eqnarray*}
Furthermore, for some constant $c>0$ and  for any constant $\alpha >0$,
\begin{eqnarray}
|\langle B(\phi) - B(\psi), \phi-\psi\rangle| &\leq & c\, |\phi-\psi|\, \|\phi - \psi \|\, \|\psi\|
\label{diffB-1}\\
&\leq & \alpha \, \|\phi - \psi\|^2 + \frac{3^3\, c^2}{2^4 \, \alpha^3}\, |\phi - \psi|^2 \, |\psi|_{L^4}^4.
\label{diffB-2}
\end{eqnarray}
\end{lemma}
\begin{proof}
Integration by parts, the boundary conditions and $div(u)
  \, = \, \nabla \cdot u
=0$ yield
\begin{eqnarray*}
 \big\langle  B_2(\phi)-B_2(\psi), \theta-\eta\big\rangle&=&\int_D \big(u.\nabla\theta)(\theta-\eta) dx
- \int_D \big(v.\nabla\eta)(\theta-\eta) dx\\
&=& - \int_D \big(u.\nabla(\theta-\eta)\big)\theta dx + \int_D \big(v.\nabla(\theta-\eta)\big)\eta dx
\end{eqnarray*}
Since $
 \langle B_2(u,w)\, ,\, w\rangle =
 \int_D \big(u.\nabla w\big) w dx =0$ for any $w\in V_2$, we deduce that
\[ \big\langle  B_2(\phi)-B_2(\psi), \theta-\eta\big\rangle =
 -\int_D  \big(u.\nabla(\theta-\eta))\eta dx + \int_D \big(v.\nabla(\theta-\eta)\big) \eta dx,\]
which completes the proof of the second identity. The proof of the first one, which is
similar and classical,  is omitted.
Finally, combining these identities with the upper estimates in Lemmas \ref{B2} and \ref{BB}
concludes the proof.
\end{proof}

For  $\phi=(u,\theta) \in V$,  define
\begin{eqnarray} \label{F}
F(\phi)=- A\phi-B(\phi)-R \phi.
\end{eqnarray}
We at first prove  crucial   monotonicity properties  of $F$.
Let $\nu \wedge\k := \min(\nu,\k)$.
\begin{lemma} \label{monotone}
 Assume that $\phi=(u,\theta)\in V$ and $\psi=(v,\eta)
\in V$; 
then  for some constant $c>0$ 
we have
\begin{equation}
\label{mono1}
  \big\langle
 F(\phi)-F(\psi), \phi-\psi  \big\rangle   + (\nu \wedge \kappa) \|\phi-\psi\|^2
\leq c\,  |\phi-\psi| \|\phi-\psi\| \|\psi\| + |\phi-\psi|^2 .
 \end{equation}
\end{lemma}
\begin{proof}
Set $U : = u-v,\; \Theta  := \th-\eta$ and $\Phi=\phi-\psi=(U, \Theta)$.
Integrating by parts 
 we deduce from Lemma \ref{diffB}
\[   \big\langle  F(\phi)-F(\psi), \Phi \big\rangle
 \; =\;  -\nu \|U\|^2-\k
 \| \Theta\|^2
 - \langle B_1(U,U),v\rangle  -   \langle B_2(\Phi),\eta\rangle
 +2(U_2,\Theta) .
\]
Thus 
\eqref{diffB-1} yields \eqref{mono1}.
\end{proof}

\bigskip

The proof of Theorem \ref{wellposeness} involves Galerkin
approximations. Let $\{\varphi_n \}_{n\geq1}$ be a complete
orthonormal basis of the Hilbert space $H$ such that $\varphi_n
\in Dom (A)$, domain of definition of the operator $A$. For any
$n\geq 1$, let  $ H_n = span(\varphi_1, \cdots, \varphi_n) \subset
Dom(A)$    and $P_n: H \to H_n$ denote the orthogonal projection
onto $H_n$.   Note that $P_n$ contracts the $H$ and $V$ norms and that its norm
as  a linear operator of $L^4(D)^3$ is bounded in $n$.
Suppose that the $H-$valued Wiener process $W$ with covariance
operator $Q$ is such that
$$
P_n Q^\frac12 = Q^\frac12 P_n, \;\; n\geq 1,
$$
which is true if $Q h=\sum_{n\geq1} \lambda_n
\varphi_n
$ with trace
$\sum_{n\geq 1}\lambda_n <\infty$. Then for $H_0 = Q^{\frac12} H$
and $ (\phi, \psi)_0 = (Q^{-\frac12}\phi, Q^{-\frac12}\psi )$ given 
$\phi, \psi \in H_0$, we see that $P_n: H_0 \to H_0 \cap H_n$ is a
contraction
 both  of   the  $H$ and  $H_0$ norms.
Let $W_n=P_n W$,  $\s_n=P_n\s$ and $\tilde{\s}_n=P_n \tilde{\s}$.

 For
  $h \in \mathcal{A}_M$,
consider the following stochastic ordinary differential equation
on the $n$-dimensional space
 $H_n$  defined  by
\begin{eqnarray} \label{ndim}
d(\phi_{n,h}^\e, \psi)=\big[ \langle F(\phi_{n,h}^\e),\psi\rangle +(\tilde{\s}_n(\phi_{n,h}^\e)h,
\psi) \big]dt + \sqrt{\e}\, (\s_n(\phi_{n,h}^\e) dW_n, \psi),
\end{eqnarray}
for $\psi=(v, \eta) \in H_n$ and $\phi_{n,h}^\e(0)=P_n \xi$.

Note that for $\psi=(v,\eta)\in V$,
  the map $\phi 
 \in H_n  \mapsto \langle (A+R)(\phi), \psi\rangle $ is
globally Lipschitz, while using Lemma
\ref{B2}  the map $\phi  =(u, \th)  \in H_n  \mapsto
\sum_{i,j=1,2} \int_D u_i \,  v_j\;\partial_i u_j \; dx  +
\sum_{i=1,2} \int_D u_i \,\eta \;\partial_i \,  \th \, dx $ is
locally Lipschitz.
 Furthermore, conditions  ({\bf A.3}) and ($\tilde{\rm\bf A}$.2) imply  that the
maps $\phi \in H_n \to \s_n(\phi)$ and $\phi \in H_n \to \tilde{\s}_n(\phi)$ are globally Lipschitz from $H_n$
to $n\times n$ matrices. Hence by a well-posedness result for
stochastic ordinary differential equations  \cite{Kunita}, there
exists a maximal solution to \eqref{ndim}, i.e., a stopping time
$\tau_{n,h}^\e\leq T$ such that \eqref{ndim} holds for $t<
\tau_{n,h}^\e$ and as $t \uparrow \tau_{n,h}^\e<T$,
$|\phi_{n,h}^\e(t)| \to \infty$.
 For every $N>0$, set
\begin{equation} \label{TauN}
\tau_N = \inf\{t:\; |\phi_{n,h}^\e(t)| \geq N \}\wedge T.
\end{equation}
Almost surely, $\phi_{n,h}^\e \in C([0, T], H_n)$ on $\{\tau_N=T \}$.
The following proposition shows that  $\tau_{n,h}^\e=T$   a.s. and  gives
estimates on $\phi_{n,h}^\e$ depending only on  the physical constants $\nu$ and $\kappa$,
$K$, $\tilde{K}$, $T$, $M$,  $\EX |\xi|^{2p}$ which are valid for all $n$ and all $\e \in [0,  \e_0]$ for some $\e_0>0$.
Its proof depends on the following version of Gronwall's lemma.
\begin{lemma} \label{lemGronwall}
Let $X$, $Y$ and  $I$ be non-decreasing, non-negative processes,  
$\varphi$ be a non-negative process and
 $Z$ be a non-negative
integrable random variable. Assume that $\int_0^T \varphi(s)\, ds \leq C$ almost surely and that there
exist positive constants $\alpha, \beta  \leq \frac{1}{2(1+Ce^C)} $,
 ${\gamma} \leq \frac{\alpha}{2(1+Ce^C)}$  and  $\tilde{C}>0$ such that
 for $0\leq t\leq T$,
\begin{eqnarray} \label{XYZ}
X(t)+ \alpha Y(t) & \leq & Z + \int_0^t \varphi(r)\, X(r)\,  dr + I(t),\;  \mbox{\rm a.s.}\\
 \label{borneI}
\EX(I(t)) &\leq & \beta\, \EX (X(t)) + {\gamma}\, \EX (Y(t))  + \tilde{C}.
\end{eqnarray}
Then if $ X \in L^\infty([0,T] \times \Omega)$,
we have for $t\in [0,T]$
\begin{equation} \label{Gronwall}
\EX \big[ X(t) + \alpha Y(t)\big]  \leq 2 (1+Ce^C) \big( \EX(Z) +\tilde{C} \big).
\end{equation}
\end{lemma}
\begin{proof}
Iterating  inequality \eqref{XYZ} and ignoring $Y$, an induction argument on $n$ yields
for $t\in [0,T]$, $n\geq 1$
\begin{align*}
 X(t) &\, \leq  \, Z  +\int_0^t \varphi(s_1)\Big[ Z
+ \int_0^{s_1} \varphi(s_2) X(s_2) ds_2 + I(s_1)\Big]
ds_1 + I(t)\\
&\leq Z + I(t)  
+ \sum_{1\leq k\leq n} \int_0^t\varphi(s_1)\int_0^{s_1} \varphi(s_2)
 \cdots \int_0^{s_{k-1}}
\varphi(s_k)\, [Z+I(s_k)]\, ds_k  \cdots ds_1 
\\
&\quad  + \int_0^t\varphi(s_1)\int_0^{s_1} \varphi(s_2) \cdots \int_0^{s_{n}}
\varphi(s_{n+1})\, X(s_{n+1})\,  ds_{n+1}  ds_n \cdots ds_1 .
\end{align*}
Recall that $X(s,\omega)$ is a.e. bounded 
 and $\int_0^T \varphi(s)\, ds \leq C$; thus
$ X(t)\leq e^C \, [Z+ I(t)]$. Using this inequality in \eqref{XYZ} and the fact that $I$
is non-decreasing, we deduce that
$ X(t)+\alpha Y(t) \leq \big[Z +I(t)\big] \,  \big( 1+Ce^C\big)$. 
Taking expected values and using \eqref{borneI}, we  conclude the proof.
\end{proof}
\begin{prop} \label{Galerkin}
There exists $\e_{0,p}:= \e_{0,p}(\nu,\kappa, K, \tilde{K},  T, M )$ such that for $0\leq \e\leq \e_{0,p}$ the
following result holds for  an integer $p\geq 1$ (with the convention $x^0=1$). 
 Let $h\in \mathcal{A}_M$ and $\xi \in L^{2p}(\Om, H)$. Then
$\tau_{n,h}=T$ a.s. and  equation \eqref{ndim} has a unique
solution with a modification  $\phi_{n,h}^\e \in C([0, T], H_n)$ and
satisfying
\begin{align} \label{Galerkin1}
\sup_n \EX \,\Big(\, & \sup_{0\leq t\leq T}|\phi_{n,h}^\e(t)|^{2p}  +
\int_0^T
\|\phi_{n,h}^\e(s)\|^2 \, |\phi_{n,h}^\e(s)|^{2(p-1)} ds \, \Big)  \nonumber \\
& \leq  C(p, \nu, \k, K, \tilde{K}, T, M) \big( \EX|\xi|^{2p} +1\big).
\end{align}
\end{prop}
\begin{proof}
  It\^o's
 formula
yields that for $t \in [0, T]$ and $\tau_N$ defined by \eqref{TauN},
\begin{align} \label{Galerkin3}
& |\phi_{n,h}^\e(t\wedge \tau_N)|^2 = |P_n \xi|^2
+ 2\sqrt{\e} \int_0^{t\wedge \tau_N}
 \big( \sigma_n(\phi_{n,h}^\e(s))  dW_n(s) , \phi_{n,h}^\e(s)\big)\\
& \qquad +2\int_0^{t\wedge \tau_N} \big\langle F(\phi_{n,h}^\e(s)), \phi_{n,h}^\e(s)\rangle ds
+ 2\int_0^{t\wedge \tau_N}
 \big( \tilde{\sigma}_n(\phi_{n,h}^\e(s)) h(s) , \phi_{n,h}^\e(s)\big)\, ds \nonumber \\
&\qquad + \e \int_0^{t\wedge \tau_N} |\sigma_n(\phi_{n,h}^\e(s))\, P_n|_{L_Q}^2\, ds.
\end{align}
Apply again It\^o's formula for $x^{p}$ when $p\geq 2$ and then  use Lemma \ref{lemBi}.
With the convention $p(p-1) x^{p-2}=0$ for $p=1$ , this yields for $t \in [0, T]$,
\begin{align} \label{estimate1}
|\phi_{n,h}^\e(t\wedge \tau_N)|^{2p} & + 2p \int_0^{t\wedge \tau_N} |\phi_{n,h}^\e(r)|^{2(p-1)} \,
 \big[ \nu  \|u_{n,h}^\e(r)\|^2+\k \|\th_{n,h}^\e(r)\|^2 \big] \, dr
 \nonumber \\
&\leq \;   |P_n\xi|^{2p}  + \sum_{1\leq j\leq 5} {T}_j(t),
\end{align}
where
\begin{eqnarray*}
{T}_1(t) &= & 4p\,  \int_0^{t\wedge \tau_N}
|(\th_{n,h}^\e(r),u^\e_{n,h,2}(r))| |\phi_{n,h}^\e(r)|^{2(p-1)} dr,   \\
{T}_2(t) &= & 2p\sqrt{\e}\;  \Big| \int_0^{t\wedge
\tau_N}
\big(\s_n(\phi_{n,h}^\e(r))\; dW_n(r),
 \phi_{n,h}^\e(r)\big )\; |\phi_{n,h}^\e(r) |^{2(p-1)} \Big| , \\
{T}_3(t) &= & 2p \,  \int_0^{t\wedge \tau_N}
|(\tilde{\s}_n(\phi_{n,h}^\e(r))\; h(r), \phi_{n,h}^\e(r))| \; |\phi_{n,h}^\e(r)|^{2(p-1)}
dr, \\
T_4(t) &= &  p\, \e \,  \int_0^{t\wedge \tau_N}
|\s_n(\phi_{n,h}^\e(r))\; P_n |^2_{L_Q} \; |\phi_{n,h}^\e(r) |^{2(p-1)}  dr,  \\
{T}_5(t) &= & 2p (p-1) \e \,  \int_0^{t\wedge \tau_N}
|\Pi_n\, \s_n^*(\phi_{n,h}^\e(r))\;  \phi_{n,h}^\e(r)|^2_{H_0}\, |\phi_{n,h}^\e(r)|^{2(p-2)} dr.
\end{eqnarray*}
The Cauchy-Schwarz inequality implies that $2 |(\th_{n,h}^\e(r),u^\e_{n,h,2}(r))|\leq
 |\phi_{n,h}^\e(r)|^2$. Hence
\begin{equation} \label{estimate2}
T_1(t)  \leq 2 p\int_0^{t\wedge \tau_N} |\phi_{n,h}^\e(r)|^{2p}\, dr.
\end{equation}
Since $h\in \mathcal{A}_M$, the Cauchy-Schwarz inequality,
 ($\tilde{\rm\bf A}$.2),  \eqref{VL4}  and the Poincar\'e inequality \eqref{Poincare}
 imply the existence of some positive
constant $c$ such that for every  $\delta_1>0$,
\begin{align} \label{estimate5}
{T}_3 (t) \leq &\; 2p  \;   \int_0^{t\wedge \tau_N} \big[\,\tilde{ K} (1+
c\, \|\phi_{n,h}^\e(r)\|^2)\big]^\frac12 \;|h(r)|_0 \;|\phi_{n,h}^\e(r)|^{2p-1} dr  \nonumber \\
\leq &\;  \delta_1\,  \int_0^{t\wedge \tau_N} \|\phi_{n,h}^\e(r)\|^2
 \, |\phi_{n,h}^\e(r)|^{2(p-1)}\, dr
+ \frac{p^2 \tilde{K}c}{\delta_1} \int_0^{t\wedge \tau_N}  |h(r)|_0^2 \,  |\phi_{n,h}^\e(r)|^{2p}\, dr
\nonumber \\
& + \delta_1\, \int_0^{t\wedge \tau_N} |\phi_{n,h}^\e(r)|^{2(p-1)} dr .
\end{align}
Using  ({\bf A.2}), we deduce that
\begin{eqnarray} \label{estimate4}
{T}_4(t) + {T}_5(t) &\leq & 2p^2 \, K\,  \e  \int_0^{t\wedge \tau_N}
 \|\phi_{n,h}^\e(r)\|^2\, |\phi_{n,h}^\e(r)|^{2(p-1)}\, dr \nonumber \\
&& + 2p^2 \, K \, \e  \int_0^{t\wedge \tau_N}  |\phi_{n,h}^\e(r)|^{2(p-1)}\, dr.
\end{eqnarray}
Finally, the Burkholder-Davies-Gundy inequality, ({\bf A.2}) and Schwarz's inequality
yield   that for $t\in[0, T]$ and $\delta_2>0$,
\begin{align} \label{estimate3}
\EX \Big(\sup_{0\leq s\leq t}|T_2(s)|\Big) \; \leq& \; 6p\sqrt{ \e} \, \EX \Big\{ \int_0^{t\wedge \tau_N}
|\phi_{n,h}^\e(r)|^{2(2p-1)} \; |\s_{n, h}(\phi_{n,h}^\e(r))\; P_n |^2_{L_Q}\;
 dr \Big\}^\frac12      \nonumber \\
\; \leq & \; \delta_2 \EX  \Big(\sup_{0\leq s\leq t\wedge\tau_N}
|\phi_{n,h}^\e(s)|^{2p}\Big)
+  \frac{9 p^2 K\e }{\delta_2}   \;\EX
 \int_0^{t\wedge \tau_N} |\phi_{n,h}^\e(r)|^{2(p-1)} dr \nonumber \\
&\quad + \frac{9 p^2 K\e }{\delta_2}
 \;\EX   \int_0^{t\wedge \tau_N}
  \|\phi_{n,h}^\e(r)\|^2 \, |\phi_{n,h}^\e(r)|^{2(p-1)} dr.
\end{align}
Consider the following property $I(i)$ for an integer $i\geq 0$:

\noindent {\bf I(i)} There exists $\e_{0,i}:= \e_{0,i}(\nu, \kappa,K, \tilde{K}, T,M)>0$ such that  for
 $0\leq \e\leq \e_{0,i}$
\[ \sup_n \EX  \int_0^{t\wedge \tau_N} |\phi_{n,h}^\e(r)|^{2i} dr \leq
C(i):= C(i,\nu, \kappa, K, \tilde{K}, T,M)<+\infty.\]
The property $I(0)$ obviously holds  with  $\e_{0,0} =1$ and $C(0)=T$.
Assume that for some integer $i$ with  $1\leq i\leq p$, the
 property I(i-1) holds;  we prove that I(i) holds.

Set $\delta_1=\frac{(\nu\wedge \kappa)\, i}{2}$, $\varphi_i(r)=
2i+\frac{i^2 \, c\tilde{ K} }{\delta_1} |h(r)|_0^2$,
$Z=
\delta_1 \, \int_0^{\tau_N} |\phi_{n,h}^\e(r)|^{2(i-1)} dr + |\xi|^{2i}$,
$X(t)=\sup_{0\leq s\leq t} |\phi_{n,h}^\e(s\wedge \tau_N)|^{2i}$,
$Y(t)=\int_0^{t\wedge \tau_N} \|\phi_{n,h}^\e(s)\|^2 \, |\phi_{n,h}^\e(s)|^{2(i-1)}\, ds$ and
$I(t)=\sup_{0\leq s\leq t} 2i\sqrt{\e}\;  \Big| \int_0^{t\wedge
\tau_N}
\big(\s_n(\phi_{n,h}^\e(r))\; dW_n(r),
 \phi_{n,h}^\e(r)\big )\; |\phi_{n,h}^\e(r) |^{2(i-1)} \Big| $.\\
Then $\int_0^T \varphi_i(s) ds\leq C_i(M):= 2iT+ \frac{i^2 c\tilde{K} }{\delta_1}M$.
Let $\alpha=i\, (\nu \wedge \kappa)$, $\beta=\delta_2=\frac{1}{2\big[1+C_i(M)e^{C_i(M)}\big]}$
 and $\tilde{C}=\frac{9 i^2 K}{\delta_2} \EX\int_0^{\tau_N} |\phi_{n,h}^\e(s)|^{2(i-1)} ds$.
Let
\[  \e_{0,i}=1\wedge \frac{\nu \wedge \kappa}{8i K}\wedge
\frac{\nu\wedge\kappa}{144 \, i\,  K\, [1+ C_i(M)e^{C_i(M)}]^2}\wedge \e_{0,i-1}.\]
Then for $0\leq \e \leq \e_{0,i}$  inequalities \eqref{estimate1}-\eqref{estimate3}
show that the assumptions of  Lemma \ref{lemGronwall} hold with
$\gamma=\frac{9 i^2 K \e}{\delta_2} \leq \alpha\beta$,
which yields I(i).

 An induction argument shows  that $I(p-1)$ holds, and hence the previous
computations with $i=p$ and  Lemma  \ref{lemGronwall} yield that for $t=T$ and    $0\leq \e\leq
\e_{0,p}$,
\[\sup_n \EX\Big( \sup_{0\leq s\leq \tau_N} |\phi_{n,h}^\e(s)|^{2p}
+ \int_0^{\tau_N} \|\phi_{n,h}^\e(s)\|^2 \, \phi_{n,h}^\e(s)|^{2(p-1)}\, ds \Big)
\leq C(p,\nu,\kappa,K,\tilde{K},  T,M).\]
As $N \to \infty$, $\tau_N \uparrow \tau_{n, h}$ and on
$\{\tau_{n, h} < T  \}$,
 $  \sup_{0\leq s \leq  t\wedge \tau_N}|\phi_{n,h}(s)| \to \infty$.
 Hence $\PX (\tau_{n, h} < T)=0$ and for almost all $\om$, for
 $N(\om)$ large enough,  $\tau_{N(\om)}(\om)=T$ and
$\phi_{n,h}(.)(\om) \in C([0, T], H_n)$.
By the Lebesgue monotone convergence theorem, we complete the
proof of the proposition.
\end{proof}
We now have the following
bound in $L^4(D)^3$.
\begin{prop} \label{L4}
  Let $h\in \mathcal{A}_M$
 and $\xi \in L^4(\Om, H)$. Let $\e_{0,2}$ be defined as in Proposition \ref{Galerkin} with $p=2$.
Then there exists a constant
  $C_{2}:= C_{2}(\nu,\kappa, K,\tilde{K},  T, M)$ 
 such that  for $0\leq \e\leq \e_{0,2}$, 
\begin{eqnarray} \label{L4bound}
 \sup_n \EX  \int_0^T |\phi_{n,h}^\e(s) |^4_{ L^4} ds
&\leq & C_{2}   (1+ \EX|\xi|^4).
\end{eqnarray}
\end{prop}
\begin{proof}
Let $f_{n,h}(t)=u_{n,h, i}(t)$ or $\th_{n,h}^\e(t)$, with
$i=1, 2$. Then \eqref{Galerkin1} with $p=2$ implies that
\begin{eqnarray*}
\sup_n \EX  \int_0^T \|f_{n,h}(s) \|^2 |f_{n,h}(s)|^2 ds \leq
 C_{2}(\nu,\kappa, K, \tilde{K},  T, M) (1+ \EX|\xi|^4).
\end{eqnarray*}
Hence by the second part of \eqref{normB1V'}, we finish the proof of
\eqref{L4bound}.
\end{proof}
\smallskip

The following result is a consequence of It\^os formula; it will be used in what followe  for various
choices of coefficients.
\begin{lemma} \label{Ito}
Let $\xi \in L^4(\Omega,H)$ be ${\mathcal F}_0$-measurable,  
$\rho' : [0,T]\times \Omega \to [0+\infty[$ be adapted such that for almost every $\omega$
the map $t\to \rho'(t,\omega)\in L^1([0,T])$ and for $t\in [0,T]$, set $\rho(t)=\int_0^t \rho'(s)\, ds$.
For $i=1,2$, let $\s_i$ satisfy assumption ({\bf A.1}), $\bar{\s}_i \in C([0,T]\times H , L^2_Q)$ and
let $\bar{\sigma}$ satisfy Assumption $\tilde{\rm\bf A}$. Let $F$ satisfy condition \eqref{mono1},
 $h_\e \in {\mathcal A}_M$
 and   $\phi_i\in L^2([0,T],V) \cap  L^\infty([0,T],H) $ a.s. and
be such that $\phi_i(0)=\xi $ and
satisfy the equation
\begin{equation} \label{phii}
d\phi_i(t)=F(\phi_i(t)) dt + \sqrt{\e} \sigma_i(t,\phi_i(t))\, dW(t)
+   \big(\bar{\s}(t,\phi_i(t))  h_\e(t)  + \bar{\s}_i(t) \big)   \, dt .
\end{equation}
Let $\Phi=\phi_1-\phi_2$ and $c_1$ and $c_2$ denote the constants in \eqref{VL4} and \eqref{Poincare}
respectively. Then for every $t\in [0,T]$,
\begin{eqnarray}\label{Itodiff}
e^{-\rho(t)}\, |\Phi(t)|^2 &\leq &   \int_0^t e^{-\rho(s)}\, \Big\{ -(\nu \wedge \kappa)\, \|\Phi(s)\|^2
+ \e \big| \s_1(s, \phi_1(s)) - \s_2(s,\phi_2(s))\big|_{L_Q^2}^2 \nonumber \\
&& + |\Phi(s)|^2\,
\Big[ -\rho'(s) +2+\frac{8\,c_1^2}{\nu\wedge \kappa} \, \|\phi_2(s)\|^2 + \frac{2 \tilde{L}\, c_1\, c_2}{\nu\wedge \kappa}\,
|h_\e(s)|_0^2\Big]\Big\}\, ds \nonumber \\
&& +\,  2\, \int_0^t e^{-\rho(s)} \big( \bar{\s}_1(s) -\bar{\s}_2(s) \, ,\, \Phi(s)\big)\, ds +I(t),
\end{eqnarray}
where  $I(t)=2 \sqrt{\e}\, \int_0^t e^{-\rho(s)} \Big( \big[ \s_1(s, \phi_1(s)) - \s_2(s,\phi_2(s))\big]\, dW(s)
\, , \,  \Phi(s)\Big)$.
\end{lemma}
\begin{proof}
It\^o's formula, \eqref{mono1}  and  condition ($\tilde{\rm\bf A}$.2)  imply that for $t\in [0,T]$,
\begin{align*}
& e^{-\rho(t)}\, |\Phi(t)|^2=  \int_0^t e^{-\rho(s)}\, \Big\{ -\rho'(s) |\Phi(s)|^2 +
 \e \big| \s_1(s, \phi_1(s)) - \s_2(s,\phi_2(s))\big|_{L_Q}^2  \\
& \quad
 + 2\left\langle F(\phi_1(s)) - F(\phi_2(s))\, ,\, \Phi(s)\right\rangle +
 2\big(  \big[ \bar{\sigma}(s,\phi_1(s)) - \bar{\sigma}(s,\phi_2(s))\big] h_\e(s)   \, ,\,
 \Phi(s) \big) \Big\}
\, ds   \\
& \quad +   \int_0^ t e^{-\rho(s)}  2\big(  \bar{\s}_1(s) -\bar{\s}_2(s)  \, ,\, \Phi(s)\big)\, ds
+I(t)  \\
& \leq   \int_0^ t e^{-\rho(s)} \Big\{ -\rho'(s)\, |\Phi(s)|^2  +
 \e \big| \s_1(s, \phi_1(s)) - \s_2(s,\phi_2(s))\big|_{L_Q}^2
-2 (\nu\wedge\kappa)\, \|\Phi(s)\|^2 \\
&\quad + 4 c_1 |\Phi(s)|\, \|\Phi(s)\|\, \|\phi_2(s)\| + 2 |\Phi (s)|^2
 + 2 \sqrt{\tilde{L}\,  c_1\, c_2 } \|\Phi(s)\|\, |h_\e(s)|_0
\, |\Phi(s)|  \Big\} ds   \\
&\quad +    \int_0^ t e^{-\rho(s)} 2\big(  \bar{\s}_1(s) -\bar{\s}_2(s) \, ,\, \Phi(s)\big)\, ds
+I(t) .
\end{align*}
The inequalities $ 4 c_1 |\Phi(s)|\, \|\Phi(s)\|\, \|\phi_2(s)\| \leq \frac{(\nu\wedge\kappa)}{2} \|\Phi(s)\|^2
+ \frac{8c_1^2}{\nu\wedge\kappa}  \|\phi_2(s)\|^2  |\Phi(s)|^2$ and
$ 2 \sqrt{\tilde{L}\,  c_1\, c_2 } \|\Phi(s)\|\, |h_\e(s)|_0 \, |\Phi(s)|
\leq \frac{(\nu\wedge\kappa)}{2} \|\Phi(s)\|^2
+  \frac{2 \tilde{L}\, c_1\, c_2}{\nu\wedge \kappa}\,
|h_\e(s)|_0^2\, |\Phi(s)|^2$ conclude the proof of \eqref{Itodiff}.
\end{proof}
We are now ready to prove the main result of this section.
\smallskip

\noindent \textbf{Proof of Theorem \ref{wellposeness}: }\\
Let $\Om_T = [0, T]\times \Om$ be endowed with the product measure
$ds \otimes d\PX$ on $ \mathcal{B} ([0, T]) \otimes \mathcal{F}$.
Let $\e_{0,2}$ be defined by Proposition \ref{Galerkin} with $p=2$ and set $\e_0:=
\e_{0,2} \wedge \frac{\nu\wedge \kappa}{2L}$.
  The proof consists of several steps.
\medskip

\noindent \textbf{Step 1: } \;
Inequalities \eqref{Galerkin1} 
and \eqref{L4bound} imply the existence of a subsequence of
$\{\phi_{n,h}^\e\}_{n\geq 0}$ (still denoted by the same notation),
 of  processes  $\phi_h^\e \in L^2(\Om_T, V) \cap L^4(\Om_T,
L^4(D)^3) \cap L^4(\Omega , L^\infty([0, T], H))$,   $F_h^\e \in L^2(\Om_T, V')$, 
 $ S^\e_h, \tilde{S}^\e_h \in L^2(\Om_T, L_Q)$,  and of random variables
 $\tilde{\phi}_h^\e(T) \in
L^2(\Om, H)$,
for which the following  properties hold:
\\
(i) $\phi_{n,h}^\e  \to \phi_h^\e$ weakly in $L^2(\Om_T, V)$, \\
(ii) $\phi_{n,h}^\e  \to \phi_h^\e$ weakly in $L^4(\Om_T,L^4(D)^3)$, \\
(iii)    $\phi_{n,h}^\e $ is weak star converging to $ \phi_h^\e$  in $L^4(\Omega, L^\infty([0, T], H))$, \\
(iv) $\phi_{n,h}^\e(T)  \to
 \tilde{\phi}_h^\e(T)
$ weakly in $L^2(\Om, H)$, \\
(v) $F(\phi_{n,h}^\e) \to F_h^\e$ weakly in $L^2(\Om_T, V')$,   \\
(vi) $\s_n(\phi_{n,h}^\e) P_n \to S_h^\e$ weakly in $L^2(\Om_T,L_Q)$,   \\
(vii)  $\tilde{\s}_n(\phi_{n,h}^\e) h \to \tilde{S}_h^\e$  weakly in $L^{\frac{4}{3}}(\Omega_T,H)$.

Indeed, (i)-(iv) are straightforward consequences of
 Propositions \ref{Galerkin} and \ref{L4}, and of
 uniqueness of the limit
of $\EX \int_0^T \phi_{n,h}^\e(t) \psi(t)dt$ for appropriate $\psi$.

Furthermore, given $\psi = (v, \eta) \in L^2(\Om_T, V)$, we have
\begin{align} \label{estimate6}
 \EX \int_0^T & \Big [\nu \langle A_1(u_{n,h}^\e(t), v(t)\rangle +\k \langle A_2(\th_{n,h}^\e(t)),
\eta(t)\rangle \big]\, dt     \nonumber \\
&= \;  -\nu \EX\int_0^T (\nabla u_{n,h}^\e(t), \nabla v(t))dt -\k\EX
\int_0^T (\nabla  \th_{n,h}^\e(t), \nabla \eta(t)) dt    \nonumber \\
& \to \;  -\nu \EX\int_0^T (\nabla u_{h}^\e(t), \nabla v(t))dt -\k\EX
\int_0^T (\nabla  \th_{h}^\e(t), \nabla \eta(t)) dt .
\end{align}
Using \eqref{Galerkin1} with $p=2$,
\eqref{normB1V'},  \eqref{inegB2},
 the Cauchy-Schwarz and Poincar\'e inequalities,
 we  deduce
\begin{align*}
 \sup_n &\;   \EX  \int_0^T \big|
 \langle  B_1(u_{n,h}^\e(t), u_{n,h}^\e(t)), v(t) \rangle
 +\langle B_2(\phi_{n,h}^\e(t))\, ,\, \eta(t)\rangle
 + (R\phi_{n,h}^\e(t), \psi(t))\big|\,  dt   \\
&\leq \; C\,  \sup_n   \EX \int_0^T \big\{
  \|u_{n,h}^\e(t)\| \; |u_{n,h}^\e(t)| \|v(t)\| +
  \|\phi_{n,h}^\e(t)\| \; |\phi_{n,h}^\e(t)| \; \|\eta(t)\|
 \\
 &\qquad\qquad \qquad  +  |\th_{n,h}^\e(t)|\; | v_2(t)| +  |u_{n,h,2}^\e(t)|\; |\eta(t))|  \big\} dt    \\
 &\leq \;  C_3( \nu, \k, K, T, M) \big(1+ E|\xi|^4\big) +
\EX \int_0^T  \|\psi(t)\|^2  dt .
\end{align*}
Hence $\{B(\phi_{n,h}^\e(t)) + R \phi_{n,h}^\e(t)\,, \,  n\geq 1 \} $ has a subsequence
converging weakly in $L^2(\Om_T, V')$. This convergence and \eqref{estimate6}  prove (v).

Since $P_n$ contracts the $|\cdot|_0$ and $|\cdot |$ norms, ({\bf A.2})
and \eqref{Galerkin1}
 imply that
\begin{eqnarray*}
\sup_n \EX \int_0^T |\s_n(\phi_{n,h}^\e(t)) P_n|^2_{L_Q} dt \leq K
\sup_n \EX    \int_0^T   (1+   \|\phi_{n,h}^\e(t)\|^2) dt\,
< \, \infty,
\end{eqnarray*}
which proves (vi).
Finally, using Assumption ($\tilde{\rm\bf A}$.1), H\"older's inequality and  \eqref{L4bound}, we deduce that
for $h\in {\mathcal A}_M$, for any $n\geq  1$,
\begin{eqnarray*}
\EX\int_0^T |\tilde{\s}_n(\phi_{n,h}^\e(s))\, h(s)|_H^{\frac{4}{3}}\,  ds  & \leq &
\EX\int_0^T \big[ \tilde{K} (1+|\phi^\e_{n,h}(s)|_{L^4}^2)\big]^{\frac{2}{3}} |h(s)|_0^{\frac{4}{3}} \, ds\\
&\leq & \tilde{K}^{\frac{4}{3}}  \Big( \EX\int_0^T |h (s)|_0^2\, ds\Big)^{\frac{2}{3}}
\Big( \EX\int_0^T [1+|\phi_{n,h}^\e(s)|_{L^4}^2]\, ds \Big)^{\frac{1}{3}}\\
&\leq & C(M,T,K,\tilde{K},\nu,\kappa).
\end{eqnarray*}
This completes the proof of (vii).
\medskip

\noindent \textbf{Step 2: }
For $\delta >0$, let $f \in H^1(-\delta, T+\delta)$ be such that
 $\|f\|_\infty = 1$,
 $f(0)=1$ and for any integer $j \geq 1$ set
 $g_j(t)=f(t)\varphi_j$,
 where $\{\varphi_j \}_{j\geq 1}$
 is
the previously chosen orthonormal basis for $H$.  It\^o's  formula
implies that for any $j\geq 1$, and for $0 \leq t \leq T$,
\begin{eqnarray} \label{equality1}
  \big( \phi_{n,h}^\e(T)\, ,\, g_j(T)\big) = \big( \phi_{n,h}^\e(0)\, ,\,  g_j(0)\big)
 +\sum_{i=1}^4  I_{n ,k}^i,
\end{eqnarray}
where
\begin{eqnarray*}
I_{n ,k}^1 & = & \int_0^T (\phi_{n,h}^\e(s), \varphi_j) f'(s) ds, \\
I_{n ,k}^2 & = & \sqrt{\e} \, \int_0^T \big( \s_n(\phi_{n,h}^\e(s)) P_n dW_n(s), g_j(s)\big), \\
I_{n ,k}^3 & = & \int_0^T \langle F(\phi_{n,h}^\e(s)), g_j(s)\rangle  ds, \\
I_{n ,k}^4 & = & \int_0^T \big( \tilde{\s}_n(\phi_{n,h}^\e(s)) h(s), g_j(s)\big) ds.
\end{eqnarray*}
Since $f' \in L^2([0, T])$ and  for every $X \in L^2(\Om)$,
  $(t, \om) \mapsto \varphi_j  X(\omega)\, f'(t)
 \in L^2(\Om, H)$,   (i) above implies that
as $n \to \infty$, $I_{n ,k}^1 \to \int_0^T (\phi_{h}^\e(s), \varphi_j)
f'(s) ds$ weakly in $L^2(\Om)$. Similarly, (v) implies that as
$n\to \infty$, $I_{n ,k}^3 \to \int_0^T \langle F_h^\e(s), g_j(s)\rangle ds$
weakly in $L^2(\Om)$, while (vii) implies that $I_{n ,k}^4 \to \int_0^T \big(
\tilde{S}_h^\e(s),
g_j(s)\big) ds$
weakly in $L^{\frac{4}{3}}(\Omega)$.
To prove the convergence of $I_{n ,k}^2$, as in \cite{Sundar}, let
$ \mathcal{P}_T$ denote the class of predictable processes in
$L^2(\Om_T, L_Q(H_0, H))$ with the inner product
\begin{eqnarray*}
 (G, J)_{\mathcal{P}_T} =\EX \int_0^T \big(G(s), J(s)\big)_{\mathcal{P}_T} ds=
 \EX \int_0^T trace (G(s)QJ(s)^*) ds.
\end{eqnarray*}
The map
 $\mathcal{T}: \mathcal{P}_T \to  L^2(\Omega)$
 defined by
$
 \mathcal{T} (G)(t) =
 \int_0^T
\big( G(s)  dW(s) , g_j(s)\big)
$
 is linear and continuous
 because of the It\^o isometry.
 Furthermore, (vi) shows that for every
 $G \in \mathcal{P}_T$, as $n\to \infty$,
 $\big(\s_n(\phi_{n,h}^\e) P_n, G\big)_{\mathcal{P}_T} \to (S_h^\e,
 G)_{\mathcal{P}_T}$ weakly in
 $L^2(\Omega)$.

Finally, as $n\to \infty$, $P_n\xi =\phi_{n,h}^\e(0) \to
\xi $ in $H$ and by (iv), $(\phi_{n,h}^\e(T), g_j(T)) \to
(\tilde{\phi}_h^\e(T), g_j(T))$ weakly in $L^2(\Om)$. Therefore,
\eqref{equality1} leads to, as $n\to \infty$,
\begin{align} \label{equality2}
  (\tilde{\phi}_{h}^\e(T), \varphi_j)\, f(T)
 &= \; \big( \xi,\varphi_j\big)
 + \int_0^T \big(\phi_{h}^\e(s), \varphi_j\big)
f'(s) ds   + \sqrt{\e} \int_0^T \big( S_h^\e(s)dW(s), g_j(s) \big)  \nonumber \\
 &\qquad  + \int_0^T \langle F_h^\e(s) , g_j(s)\rangle ds +\int_0^T \big(
 \tilde{S}_h^\e(s) , g_j(s)\big)
ds.
\end{align}

For $\delta >0$, $k>\frac{1}{\delta}$, $t \in [0, T]$, let $f_k \in H^1(-\delta, T+\delta)$ be such
that
 $\|f_k\|_\infty =1$,
  $f_k=1$ on $(-\delta, t-\frac{1}{k})$ and $f_k=0$ on
$\big(t, T+\delta\big)$.
 Then $f_k \to 1_{(-\delta, t)}$ in
$L^2$, and $f'_k \to -\delta_t$ in the sense of distributions.
Hence as $k\to \infty$, \eqref{equality2} written with $f:=f_k$
yields
 \begin{eqnarray*}
0&=& \big(\xi ,\varphi_j\big)
- \big({\phi}_{h}^\e(t) ,\varphi_j \big)
  + \sqrt{\e}  \Big( \int_0^t S_h^\e(s)dW(s), \varphi_j\Big)
  \nonumber \\
&&
 +  \Big\langle \int_0^t   F_h^\e(s)\, ds \, , \,  \varphi_j\Big\rangle
+\Big( \int_0^t \tilde{S}_h^\e(s)\, ds\, , \, \varphi_j\Big) .
\end{eqnarray*}
Note that $j$ is arbitrary and $\EX \int_0^T |S_h^\e(s)|^2_{L_Q} ds <
\infty$;  we deduce that for $0\leq t \leq T$,
\begin{eqnarray} \label{equality3}
 \phi_{h}^\e (t)  =   \xi
 + \sqrt{\e} \int_0^t  S_h^\e(s)dW(s)  + \int_0^t F_h^\e(s)  ds +\int_0^t  \tilde{S}_h^\e(s)ds \in H.
\end{eqnarray}
Indeed, $\int_0^t F^\e_h(s)  ds$, as a  linear combination of $H-$valued
terms, also belongs to $H$.
Moreover, let $f=1_{(-\delta,
T+\delta)}$. Using \eqref{equality2} again, we obtain
\begin{eqnarray*}
 \tilde{\phi}_{h}^\e (T)
  =   \xi
 + \sqrt{\e} \int_0^T  S_h^\e(s)dW(s)  + \int_0^T F_h^\e(s)  ds +\int_0^T \tilde{S}_h^\e(s)ds.
\end{eqnarray*}
This equation and \eqref{equality3} yield that
 $\tilde{\phi}_h^\e(T) = \phi_h^\e(T)$ a.s.
\medskip

\noindent \textbf{Step 3: }
In \eqref{equality3} we still have to prove that  $ds \otimes d\PX$ a.s. on $\Om_T$,  one has
\begin{eqnarray*}
 S^\e_h(s)=\s(\phi_h^\e(s)), \;     F_h^\e(s)=F(\phi_h^\e(s))\;\mbox{\rm and }\;
 \tilde{S}^\e_h(s)=\tilde{\s}(\phi_h^\e(s))\;h(s) .
\end{eqnarray*}
Let
\begin{eqnarray*}
   {\mathcal X}&:=& \{ \psi\in  L^4(\Om_T, L^4(D)^3)\cap L^4\big( \Omega,
 L^\infty([0,T],H) \big) \cap L^2(\Omega_T,V)\; : \\
&&  \int_0^T \big( \|\psi(t)\|^2 + \|\phi_h^\e(t)\|^2\big)   |\psi(t)-\phi^\e_h(t)|^2
\, dt <+\infty\; a.s.\}.
\end{eqnarray*}
Then (i)-(iii)  yield  $\phi_h^\e\in  {\mathcal X}$ and since $\|u\|\leq C(m) |u|$ for every $u\in H_m$, using
 \eqref{normB1V'} and the fact that $\phi_h^\e\in L^2(\Omega_T, V)$,  we deduce that
for any $m\geq 1$,  $L^\infty(\Omega_T,H_m)
\subset  {\mathcal X}$.  Let  $\psi=(v, \eta)\in   L^\infty(\Omega_T,H_m)$.
 For every $t \in [0, T]$,
 if $a\wedge b := \inf(a,b)$ and $c_1$ is the constant in \eqref{VL4},
set
\begin{eqnarray} \label{r}
 r (t)  = \int_0^t  \Big[\,
 2+  \frac{8 c_1^2}{\nu \wedge\k} \, \|\psi(s)\|^2
 +  \frac{2 c_1 c_2 \tilde{L} }{\nu\wedge \kappa}  |h(s)|_0^2   \Big]\, ds.
\end{eqnarray}
Then $r(t) < \infty$ for all $t \in [0, T]$ and Fatou's lemma implies
\[ \EX\big( |\phi_h^\e(T)|^2\, e^{-r(T)}\big) \leq   \EX\big( \liminf_n |\phi_{n,h}^\e(T)|^2\, e^{-r(T)}\big)
\leq \liminf_n \EX\big(|\phi_{n,h}^\e(T)|^2\, e^{-r(T)}\big).\]
Apply  It\^o's formula  to \eqref{equality3} and  \eqref{ndim},   and for $\phi=\phi_h^\e$ or
$\phi=\phi_{n,h}^\e$, let $\phi = \psi + (\phi-\psi)$. After simplification, this yields
\begin{align} \label{r1}
\EX &|\xi|^2 + \EX\int_0^T \!\!  e^{-r(s)}\,  \big[ -r'(s)\big\{ \big|\phi_h^\e(s)-\psi(s)\big|^2
+   2\big( \phi_h^\e(s)-\psi(s)\, ,\, \psi(s)\big)\}
 + 2\langle F_h^\e(s),\phi_h^\e(s)\rangle \nonumber \\
& 
+\e |S^\e_h(s)|^2_{L^2_Q} + 2\big(\tilde{S}^\e_h(s)\, ,\, \phi_h^\e(s)\big)\big]\, ds
\leq \liminf_n \big( \EX|P_n(\xi)|^2 + X_n\big) ,
\end{align}
where
\begin{align*}
& X_n=\EX\int_0^T e^{-r(s)}\big[  -r'(s) \big\{ \big|\phi_{n,h}^\e(s)-\psi(s)\big|^2
+ 2\big( \phi_{n,h}^\e(s)-\psi(s)\, ,\, \psi(s)\big) \big\}  \\
& 
  + 2\langle F(\phi_{n,h}^\e(s)),\phi_{n,h}^\e(s)\rangle
+\e | \sigma_n (\phi_{n,h}^\e(s))P_n|^2_{L^2_Q}
+ 2\big(\tilde{\sigma}(\phi_{n,h}^\e(s)) h(s)\, ,\, \phi_{n,h}^\e(s)\big)\big]\, ds .
\end{align*}
Set $a\vee b: =\max(a, b)$. Inequalities \eqref{mono1}, ({\bf A.3}),
 ($\tilde{\rm\bf A}$.2),  \eqref{r},
the Poincar\'e  and Schwarz inequalities
imply that for $0\leq \e\leq \e_0 \leq \frac{\nu\wedge \kappa}{2L}$,
\begin{align} \label{r2}
Y_n &\; : = \; \EX \int_0^T e^{-r(s)}\big [-r'(s)
|\phi_{n,h}^\e(s)-\psi(s)|^2    \nonumber \\ &
 +  \Big[ 2 \langle F(\phi_{n,h}^\e(s))-F(\psi(s)), \phi_{n,h}^\e(s)-\psi(s)\rangle
\; +\e |\s_{n}(\phi_{n,h}^\e(s))\; P_n - \s_n(\psi(s)) \; P_n|^2_{L_Q}
  \nonumber \\
&\;  + 2 \big(\big\{\tilde{\s}_{n}(\phi_{n,h}^\e(s))- \tilde{\s}_n(\psi(s))\big\}\; h(s),
\phi_{n,h}^\e(s)-\psi(s)\big) \Big] ds     \nonumber \\
\leq &\;  \EX \int_0^T e^{-r(s)} \; |\phi_{n,h}^\e(s)-\psi(s)|^2
\Big\{-r'(s)  
+ 2  + \frac{8 c_1^2}{\nu \wedge \k}
\|\psi(s)\|^2    
 +\frac{2 c_1 c_2 \tilde{L}}{\nu\wedge \kappa} |h(s)|_0^2 \Big\} ds \nonumber\\
\leq & \; 0.
\end{align}
Furthermore, $X_n=Y_n + \sum_{i=1}^2 Z_n^i$, with
\begin{align*}
& Z_n^1 \, = \,   \EX  \int_0^T \!\!  e^{-r(s)} \Big[ -2  r'(s) (\phi_{n,h}^\e(s)-\psi(s), \psi(s))
 + 2\langle F(\phi_{n,h}^\e(s)),\psi(s) \rangle
\\
\; &+ 2 \langle F(\psi(s)),\phi_{n,h}^\e(s)\rangle
 -  2\langle F(\psi(s)), \psi(s)\rangle   
+2\e \big( \sigma_n(\phi_{n,h}^\e(s))P_n \, , \, \sigma(\psi(s)\big)_{L_Q}
 \\ &
 \;  +  2\big( \tilde{\s}_{n}(\phi_{n,h}^\e(s)\big)h(s), \psi(s)\big) 
 +  2 \big( \tilde{\s} (\psi(s))\; h(s),
\phi_{n,h}^\e(s)) -2 (P_n\tilde{\s}(\psi(s)) h(s), \psi(s) \big) \, \Big]\,  ds , \\
& Z_n^2  \, =  \EX  \int_0^T \!\! e^{-r(s)} \Big[2\e  \big(\s_n(\phi_{n,h}^\e(s)) P_n,
[\s(\psi(s)) P_n-\s(\psi(s))]\big)_{L_Q}  -  \e |P_n \s(\psi(s)) P_n|^2_{L_Q}\Big]  ds.
\end{align*}
The weak convergence properties (i)-(vii) imply that, as $n\to
\infty$,  $Z_n^1 \to Z^1$ where
\begin{align} \label{r3}
 Z^1 \,
 =   \, \EX & \int_0^T \!  e^{-r(s)} \big[ -2  r'(s) \big(\phi_{h}^\e(s) - \psi(s), \psi(s)\big)
+  2\langle F_h^\e(s),\psi(s)\rangle
+ 2 \langle F(\psi(s)),\phi_{h}^\e(s)\rangle    \nonumber  \\
& - 2\langle F(\psi(s)), \psi(s)\rangle 
 +2\e\, \big( S_h^\e(s)\, ,\, \sigma(\psi(s))\big)_{L_Q} +  2 (\tilde{S}_h^\e(s), \psi(s))    \nonumber \\
& + 2 \big(\tilde{\s} (\psi(s))\; h(s), \phi_{h}^\e(s)\big)
 -2 \big(\tilde{\s}(\psi(s)) h(s), \psi(s) \big)\, \big] ds.
\end{align}
Now we study  $(Z_n^2)$; when $n\to \infty$,   $|\s(\psi(s))
(P_n-Id_{H_0})|_{L_Q}  \to 0$ a.s., and by ({\bf A.2}),
$$
\EX \int_0^T e^{-r(s)} \sup_n |\s(\psi(s))
(P_n-Id_{H_0})|_{L_Q}^2ds < \infty.
$$
Hence the Lebesgue dominated convergence theorem implies that, as
$n\to \infty$,
\begin{eqnarray*}
\EX \int_0^T e^{-r(s)}   |\s(\psi(s)) (P_n-Id_{H_0})|_{L_Q}^2ds
\to 0.
\end{eqnarray*}
Since
$ \sup_n \EX  \int_0^T e^{-r(s)} |\s_n(\phi_{n,h}^\e(s)) P_n
|^2_{L_Q} ds < \infty,$
 by the Cauchy-Schwarz inequality, we see that $Z_n^2 \to - \e \EX \int_0^T e^{-r(s)} |\sigma(\psi(s))|_{L_Q}^2 ds$.

Thus, \eqref{r1}-\eqref{r3} imply that for any $m\geq 1$ and any $\psi \in
L^\infty(\Om_T, H_m)$,
\begin{align} \label{r7}
&    \EX  \int_0^T e^{-r(s)}\Big\{ - r'(s) |\phi_h^\e (s)-\psi(s)|^2
 + 2\langle F_h^\e(s)-F(\psi(s)), \phi_h^\e(s)-\psi(s)\rangle \nonumber \\
&\;   + \e |S^\e_h(s)-\s(\psi(s))|^2_{L_Q}
+ 2\Big( \tilde{S}_h^\e(s)-\tilde{\s}(\psi(s)) h(s)\, ,\,  \phi_h^\e(s)-\psi(s)\Big) \Big\} ds
\leq  0.
\end{align}
By a density argument, this inequality extends to all $\psi \in {\mathcal X}$.
Taking $\psi=\phi_h^\e \in {\mathcal X} $,
we conclude that $S_h^\e(s) =\s(\phi_h^\e (s))$ $ ds\otimes d \PX\; a.e. $
For a real number $\lambda$, $\tilde{\psi}=(v, \eta) \in L^\infty(\Omega_T,H_m)$ for some $m$,
 set  $\psi_\lambda =\phi_h^\e -\lambda
\tilde{\psi}\in {\mathcal X}$. Thus applying  \eqref{r7}  to $\psi_\lambda$ and
neglecting $\e| \s(\phi^\e_h(s)) - \s(\psi_\lambda(s))|^2_{L_Q}$, we obtain
\begin{align} \label{r8}
 \EX  \int_0^T e^{-r(s)} & \Big[ -\lambda^2  r'(s) |\tilde{\psi}(s)|^2
 + 2 \lambda  \Big\{\langle F_h^\e(s)-F(\psi_\lambda(s)),
\tilde{\psi}(s)\rangle  \nonumber \\
&\quad
+ \Big(\tilde{S}_h^\e(s)-\tilde{\s}(\psi_\lambda(s)) h(s), \tilde{\psi}(s)\Big) \Big\}
 \Big]\, ds \leq  0.
\end{align}
Using ($\tilde{\rm\bf A}$.2), \eqref{VL4} and \eqref{Poincare}, we have for almost every $(s, \om) \in \Om_T$ as
$\lambda \to 0$,
$$
\big|\big( \big[\tilde{\s}(\psi_\lambda(s))- \tilde{\s}(\phi_h^\e(s))]\,  h(s)\, ,\,
  \tilde{\psi}(s)\big)\big| \leq
\sqrt{\tilde{L\, c_1\, c_2}}\, \lambda \,  \|\tilde{\psi}(s)\|\, |h(s)|_{0}  \,  |\tilde{\psi}(s)|
 \to 0.
$$
Furthermore, ($\tilde{\rm\bf A}$.1) \eqref{VL4} and \eqref{Poincare} imply  that for some constant $c>0$, 
\begin{align*}
 \EX  \int_0^T& \sup_{|\lambda| \leq 1}
\Big|\Big( \tilde{\s}(\psi_\lambda(s)) \; h(s), \tilde{\psi}(s)\Big)\Big|  ds   \\
 & \leq  \; \sqrt{\tilde{K} c} \, \EX \int_0^T \big(1+2 \| \phi_h^\e(s)\|^2+
2 \|\tilde{\psi}(s))\|^2 \big)^\frac12  \;  |h(s)|_0
 \; | \tilde{\psi}(s)  |  ds \\
&\leq \; c\, \tilde{K}  M + \EX\int_0^T \Big[ \big\{ 1
 +  2  \|  \phi_h^\e(s)\|^2+
2 \|\tilde{\psi}(s)\|^2 \big\} \, |\tilde{\psi}(s)|^2\Big]   ds < \infty.
\end{align*}
Hence, the Lebesgue dominated convergence theorem yields, as
$\lambda \to 0$,
\begin{eqnarray*}
 \EX \int_0^T \! \Big( \big\{ \tilde{S}^\e_h(s)-\tilde{\s}(\psi_\lambda(s)) \big\} h(s), \tilde{\psi}(s)\Big)  ds
 \to \EX \int_0^T \! \Big(\big\{ \tilde{S}_h^\e(s)-\tilde{\s}(\phi_h(s))\big\} h(s), \tilde{\psi}(s)\Big)
 ds .
\end{eqnarray*}
Furthermore,  \eqref{mono1}   yields  for $\lambda \neq 0$
\[  \big|\langle  F(\psi_\lambda 
 (s)) -F(\phi_h^\e(s)), \tilde{\psi}(s)\rangle\big| \leq
\lambda^2 \Big[( \nu \wedge \kappa)\,   \|\tilde{\psi}(s)\|^2
+ 2 c_1  \|\tilde{\psi}(s)\|^2 |\tilde{\psi}(s)| +|\tilde{\psi}(s)|^2\Big].
\]
 Using again the dominated convergence theorem, we deduce  as $\lambda \to 0$,
\begin{eqnarray*}
 \EX \int_0^T  \langle F_h^\e(s)- F(\psi_\lambda(s)), \tilde{\psi}(s)\rangle
ds  \to \EX \int_0^T  \langle F_h^\e(s)- F(\phi_h^\e(s)), \tilde{\psi}(s)\rangle  ds.
\end{eqnarray*}
Thus, dividing \eqref{r8} by $\lambda>0$ and letting $\lambda \to
0$ we obtain that for every $m$ and
 $\tilde{\psi} \in  L^\infty(\Omega_T,H_m)$,
\begin{eqnarray*}
 \EX \int_0^T \Big[ \langle F_h^\e(s)- F(\phi_h^\e(s)), \tilde{\psi}(s)\rangle
  + \big( \big\{ \tilde{S}_h^\e(s)-\tilde{\s}(\phi_h^\e(s))\big\} h(s), \tilde{\psi}(s)\big) \Big] ds \leq 0,
\end{eqnarray*}
while a similar calculation for $\lambda <0$ yields the opposite
inequality. Therefore for almost every $(s, \om) \in \Om_T$, for
every $\tilde{\psi}$ in a dense subset of $L^2(\Om_T, V)$,
\begin{eqnarray} \label{r9}
 \EX \int_0^T \Big[  \big\langle
   F_h^\e(s)- F(\phi_h^\e(s))\, , \, \tilde{\psi}(s)\big\rangle
  + \big( \big\{\tilde{S}_h^\e(s)-\tilde{\s}(\phi_h^\e(s))\big\} h(s),
\tilde{\psi}(s)\big) \Big] \,  ds =
 0.
\end{eqnarray}
Hence a.e. for $t \in [0, T]$, \eqref{equality3} can be rewritten
as
\begin{eqnarray} \label{r9bis}
 \phi_{h}^\e (t)  =   \xi
 +\sqrt{\e}  \int_0^t  \s(\phi_h^\e(s))dW_s  +
 \int_0^t \big[F(\phi_h^\e(s)) + \tilde{\s}(\phi_h^\e(s)) h(s) \big]ds.
\end{eqnarray}
Furthermore, (i), (iv) and \eqref{Galerkin1} for $p=2$ imply that
\begin{eqnarray} \label{boundgeneral1}
E\Big(  \int_0^T \|\phi_h^\e(t)\|^2\, dt \Big)& \leq& \sup_n \EX
\int_0^T \|\phi_{n,h}^\e(t)\|^2 dt \leq C
 \big( 1+E|\xi|^4\big),\\
 \label{boundgeneral2}
E\Big( \sup_{0\leq t\leq T} |\phi_h^\e(t)|^4 \big) &\leq & \sup_n \EX
\Big(\sup_{0\leq t\leq T} |\phi_{n,h}^\e(t)|^4 \Big)  \leq    C\,
\big( 1+E|\xi|^4\big).
\end{eqnarray}
Since $|x|^2\leq 1 \vee |x|^4$ for any $x\in \RR$, this completes the proof of  \eqref{boundgeneral}.
\medskip

\noindent \textbf{Step 4: }
To complete the proof of Theorem \ref{wellposeness}, we  show
that $\phi_h^\e$ has a $C([0, T], H)$-valued modification and that
the solution to \eqref{r9bis} is unique in $X:=C([0, T], H) \cap L^2([0,
T], V)$. Note that \eqref{boundgeneral} implies that if
$\tilde{\tau}_N =\inf \{t \geq 0: |\phi_h^\e(t) | \geq N \} \wedge T$
for $N>0$, $\PX (\tilde{\tau}_N <T) \leq C N^{-2}$. The
Borel-Cantelli lemma yields $\tilde{\tau}_N \to T$ a.s. when $N \to
\infty$.

 We at first prove uniqueness. Let $\psi=(v, \eta) \in
 X   
$ be another solution to
 \eqref{r9bis}.
 Then if $\bar{\tau}_N =\inf \{t \geq 0:
|\psi(t) | \geq N \} \wedge T$ for $N>0$, since $|\psi(.)|$ is
a.s. bounded on $[0,T]$, as $N\to \infty$, we have  $\bar{\tau}_N \to T$ a.s.  and
 hence $ 
{\tau}_N  =  \tilde{\tau}_N  \wedge  \bar{\tau}_N  \to T,$ \;   a.s. 

 Let $\phi_h^\e=(u_h^\e,\theta_h^\e)$, $\Phi=\phi_h^\e -\psi $,
 and $a=\frac{8 c_1^2}{\nu \wedge \kappa}$, where $c_1$ is the constant
defined in \eqref{VL4}.
Set $\rho'(t):=a\|\psi(t)\|^2$, 
 $h_\e:=h$,
$\sigma_1=\sigma_2=\sigma$, $\bar{\s}=\tilde{\s}$,
$\bar{\s}_1=\bar{\s}_2=0$. Then  $\phi_1=\phi_h^\e$ and $\phi_2=\psi$ satisfy \eqref{phii}.
Set
\[
{\mathcal I}(t)=\sup_{\tau\leq t} 2\, \sqrt{\e}\, \int_0^\tau e^{-a\int_0^s \|\psi(r)\|^2 dr}
\Big( \big[\s(\phi_h^\e(s)) - \s(\psi(s))\big]\, dW(s)\, ,\, \Phi(s)\Big),
\]
Then using Lemma \ref{Ito} and  condition ({\bf A.3}) yields for $0\leq \e\leq \e_0\leq \frac{\nu\wedge\kappa}{2L}$
\begin{eqnarray*}
\zeta(t):&= &e^{-\rho(t\wedge \tau_N)}
 |\Phi(t\wedge \tau_N)|^2 \\
& \leq &  {\mathcal I}(t\wedge \tau_N)+  \int_0^{t\wedge \tau_N} e^{-\rho(s)}  \Big\{
 \big[\e L -\nu\wedge \kappa\big]\, \|\Phi(s)\|^2
\\
&& + |\Phi(s)|^2 \big[ -a\|\psi(s)\|^2 +2 +\frac{8 c_1^2}{\nu\wedge\kappa}\|\psi(s)\|^2 +
\frac{2\tilde{L}c_1c_2}{\nu\wedge\kappa} |h(s)|_0^2\big] \Big\} ds .
\end{eqnarray*}
Thus
\[\zeta(t)+\frac{\nu\wedge \kappa}{2} Y(t) \leq \int_0^t \Big( \frac{ 2\tilde{L}\, c_1\, c_2}{\nu\wedge \kappa}\;
|h(s\wedge  \tau_n)|_0^2+2\Big)
\zeta(s)\, ds + {\mathcal I}(t\wedge \tau_n), \]
where  $Y(t)= \int_0^{t\wedge \tau_n} e^{-\rho(s)} \, \|\Phi(s)\|^2\, ds$.
Burkholder's inequality and Assumption ({\bf A.3}) imply that for all $\beta >0$ and $\e\in [0,\e_0]$,
\[ 
\EX {\mathcal I}(t\wedge\tau_n)\leq  6\sqrt{\e_0} \EX \Big( \int_0^{t\wedge\tau_N}\!\!\!
e^{-2\rho(s)}
L \|\Phi(s)\|^2\, |\Phi(s)|^2\, ds \Big)^{\frac{1}{2}}
\leq \beta \EX \sup_{0\leq s\leq t}\zeta(s) + \frac{9L\e_0}{\beta} \EX Y(t).
\]  
Since  
$\int_0^T \big( \frac{2\tilde{L} c_1 c_2}{\nu\wedge\kappa} |h(s\wedge \tau_N)|^2_0 +2\big) \, ds \leq
\frac{2 \,M\, \tilde{L}\, c_1 c_2}{\nu\wedge \kappa} +2T:=C$, Lemma \ref{lemGronwall} implies that for
$\beta= \big( 2[1+Ce^C])^{-1}$ and  $\e_0\, L$ small
enough to have  $\frac{9 \,\e_0\, L}{\beta} \leq \frac{\nu\wedge \kappa}{2}\,  \beta$, one has
\begin{equation}
 \; \EX\; \sup_{0\leq s\leq T} e^{- a\int_0^{s\wedge \tau_N}
 \|\psi(r)\|^2dr} \;  |\Phi(s \wedge \tau_N)|^2  =0.
 \end{equation}
Since
 $\lim_{N\to \infty} \tau_N =T$
 a.s., we thus deduce  $|\Phi(s, \om)| =0$ a.s. on $\Omega_T$. Thus
if $\phi_h^\e$ is in $C([0, T], H)$, we conclude that
$\phi_h^\e(t)=\psi(t)$, a.s., for every $t \in [0, T]$.
\smallskip

Finally, set
\begin{equation}
 \tilde{\rho}'(t) =  \frac{8 c_1^2}{  \nu \wedge \k}
 \;  \|\phi_h^\e(s)\|^2    +2+   \frac{2 \tilde{ L}c_1c_2}{\nu\wedge \kappa} \; |h(s)|_0^2,
 \end{equation}
 let $h_\e:=h$,  $\sigma_1=P_n\s P_n$, $\sigma_2=\sigma$, $\bar{\s}_1=0$,  
$ \bar{\s}_2(s) =\big[ \tilde{\s}(\phi_h^\e(s))
- P_n \tilde{\s}(\phi_h^\e(s)) \big] \, h(s)$
and  $\bar{\s} = P_n \tilde{\s}$.
 Then $\tilde{\rho}(t)=\int_0^t \tilde{\rho}'(s)\, ds <+\infty$ a.s.
Then $\phi_1=\phi_{n,h}^\e$ and $\phi_2=\phi_h^\e$ satisfy \eqref{phii}.
 Set $\Phi_{n,h}^\e= \phi_{n,h}^\e - \phi_h^\e$
and let  $0\leq \e \leq \e_0  \leq \frac{\nu\wedge \kappa}{4L}$.
By Lemma \ref{Ito} and condition ({\bf A.3}), we deduce that for every $t\in [0,T]$,
\begin{align*}
 &\EX  \big( e^{-\tilde{\rho}(t)} |\Phi_{n,h}^\e(t)|^2 \big) \leq
\EX\int_0^t\!  e^{-\tilde{\rho}(s)} \Big\{ \big[ 2\e L - (\nu\wedge \kappa)\big] \|\Phi_{n,h}^\e\|^2
+2 \e |P_n\s (\phi^\e_h(s))P_n -\s(\phi_h^\e(s))|_{L_Q}^2 \\
&\quad  \quad + |\Phi_{n,h}^\e(s)|^2 \big[ -\tilde{\rho}'(s) +2
+  \frac{8\, c_1^2}{\nu\wedge\kappa} \|\phi_h^\e(s)\|^2 + \frac{2\tilde{L}\, c_1\, c_2}{\nu\wedge\kappa}
|h(s)|_0^2  \big]\Big\}\, ds \\
&\quad \quad + \EX\int_0^t e^{-\tilde{\rho}(s) }\, 
2 \,  |\Phi_{n,h}^\e(s)|\, |P_n\tilde{\s}(\phi_h^\e(s)) -  \tilde{\s} 
(\phi_h^\e(s)|_{L_Q}\, |h(s)|_0\, ds  \\
& \leq \mathcal{R}(t,n) - \frac{\nu\wedge \kappa}{2} \EX \int_0^t e^{-\tilde{\rho}(s)} \|\Phi_{n,h}^\e(s)\|^2
\, ds ,
\end{align*}
where
\[ {\mathcal R}(t,n)=  \EX\int_0^t \! 
\big[ 
2 \e |P_n\s(\phi_h^\e(s))P_n - \s(\phi_h^\e(s))|^2_{L_Q}  
 +   |P_n \tilde{\s}(\phi_h^\e(s)) - \tilde{\s}(\phi_h^\e(s))|_{L_Q}^2
\big] ds , \]
and the last inequality follows from Schwarz's inequality and the definition of $\tilde{\rho}$.

Furthermore, for almost every $(s,\omega)$,  one has
$|P_n\s(\phi_h^\e(s))P_n - \s(\phi_h^\e(s))|_{L_Q} \to 0$ and
$|P_n\tilde{\s}(\phi_h^\e(s))-\tilde{\s}(\phi_h^\e(s))|_{L_Q} \to 0$ as $n\to \infty$.
Thus the dominated convergence theorem shows that $\lim_n \sup_t {\mathcal R}(t,n)\to 0$, and
thus that  $\lim_{n\to \infty} I(n)=0$, where
$$
I(n)=  \sup_{0\leq t\leq T} \EX \big(e^{-\tilde{\rho}(t)}
|\Phi_{n,h}^\e(t))|^2\big) +\EX \int_0^T e^{-\tilde{\rho}(s)} \|\Phi_{n,h}^\e(s)\|^2\, ds  .
$$
Using again  Lemma \ref{Ito} and the Burkholder-Davies-Gundy inequality,
 a similar computation  yields that for $0\leq \e\leq \e_0\leq \frac{\nu\wedge \kappa}{4L}$:
\begin{align*}
 &\EX  \Big(\sup_{0\leq t\leq T} e^{-\tilde{\rho}(t)} |\Phi_{n,h}^\e(t)|^2 \Big) \leq
\frac{1}{2}   \EX  \Big(\sup_{0\leq t\leq T} e^{-\tilde{\rho}(t)} |\Phi_{n,h}^\e(t)|^2 \Big)
\\&\quad
+ 18\e \EX\int_0^T e^{-\tilde{\rho}(s)}\,  |\s_n(\phi_{n,h}^\e(s))P_n - \s(\phi_h^\e(s))|^2_{L_Q} ds
\\&\quad +   \EX 
\int_0^T 
\Big[ 2 \e |P_n\s (\phi^\e_h(s))P_n -\s(\phi_h^\e(s))|_{L_Q}^2 
 + |P_n\tilde{\s}(\phi_h^\e(s)) -  \tilde{\s } (\phi_h^\e(s)|_{L_Q}^2\big]\, ds\\
&
 \leq  C \big[I(n) +  {\mathcal R}(T, n)\big] .
\end{align*}
Therefore, $\phi_{n,h}^\e$ has a subsequence converging a.s.
uniformly to $\phi_h^\e$ in $H$. Since  $\phi_{n,h}^\e \in C([0, T],
H)$, we conclude that $\phi_h^\e$ has a modification in $C([0, T], H)$.
\hfill $\Box$

\section{Large deviations}  \label{s4}
  We consider large deviations via a    weak convergence approach
  \cite{BD00, BD07}, based on variational representations of
  infinite-dimensional    Wiener processes.
The solution to the stochastic B\'enard model \eqref{Benard} is denoted as 
$\phi^\e = {\mathcal G}^\e(\sqrt{\e} W)$
 for a  Borel measurable function ${\mathcal G}^\e: C([0, T], H) \to X.$ The  space
$X=C([0, T]; H) \cap L^2((0, T); V)$ endowed with the metric
associated with the norm defined in \eqref{norm} is Polish. Let
$\mathcal{B}(X)$ denote its  Borel $\s-$field. We recall some
classical   definitions.
\begin{defn} 
   The random family
$\{\phi^\e \}$ is said to satisfy a large deviation principle on
$X$  with the good rate function $I$ if the following conditions hold:\\
\indent \textbf{$I$ is a good rate function.} The  function $I: X \to [0, \infty]$ is
such that for each $M\in [0,\infty[$ the level set $\{\phi \in X: I(\phi) \leq M
\}$ is a    compact subset of $X$. \\
 For $A\in \mathcal{B}(X)$, set $I(A)=\inf_{\phi \in A} I(\phi)$.\\
\indent  \textbf{Large deviation upper bound.} For each closed subset
$F$ of $X$:
$$
\lim\sup_{\e\to 0}\; \e \log \PX(\phi^\e \in F) \leq -I(F).
$$
\indent  \textbf{Large deviation lower bound.} For each open subset $G$
of $X$:
$$
\lim\inf_{\e\to 0}\; \e \log \PX(\phi^\e \in G) \geq -I(G).
$$
\end{defn}
To establish the large deviation principle, we need to strengthen the hypothesis on the
growth condition and Lipschitz property of $\sigma$ (and $\tilde{\s}$)  as follows:\\
\noindent {\bf Assumption A Bis} There exist   positive constants $K$ and $L$ such that \\
  ({\bf A.4})  $|\sigma(t,\phi)|_{L_Q}^2 \leq K\, (1+|\phi|^2)$, $\forall  t\in [0,T]$, $\forall  \phi\in V$.\\
({\bf A.5}) $|\sigma(t,\phi)-\sigma(t,\psi)|^2_{L_Q} \leq L\, |\phi-\psi|^2$, $\forall  t\in [0,T]$,
$\forall  \phi,\psi \in V$.\\
\smallskip

 Note that due to the continuous embedding $V
\hookrightarrow H$,  assumptions ({\bf A.4-A.5}) imply ({\bf A.2-A.3}) as well as 
($\tilde{\rm\bf A}$.1-$\tilde{\rm\bf A}$.2). Thus
the 
conclusions of Theorem \ref{wellposeness}  hold  if $\tilde{\s}=\s$ satisfy  
assumptions ({\bf A.4-A.5}).

The proof of the large deviation principle will use the following
technical lemma which studies time increments of the solution to
the stochastic control equation. For any integer $k=0, \cdots,
2^n-1$, and $s\in [k T 2^{-n}, (k+1) T 2^{-n}[$, set
$\underline{s}_n= kT2^{-n}$ and $\bar{s}_n=(k+1)T 2^n$. Given
$N>0$,
 $h\in {\mathcal A}_M$,
$\e\geq 0$ small enough, let  $\phi_h^\e$ denote the solution to \eqref{Benardgene} given by
Theorem \ref{wellposeness}, and for  $t\in [0,T]$,  let
\[ G_N(t)=\Big\{ \omega \, :\, \Big (\sup_{0\leq s\leq t}  |\phi_h^\e(s)(\omega)|^2 \Big)\vee
\Big(  \int_0^t \|\phi_h^\e(s)(\omega)\|^2 ds \Big) \leq N\Big\}.\]
\begin{lemma} \label{timeincrement}
Let $M,N>0$, $\sigma$ and $\tilde{\sigma}$ satisfy  Assumptions ({\bf A.1}),({\bf A.4}) and ({\bf A.5}),
 $\xi\in L^4(\Omega,H)$  be ${\mathcal F}_0$-measurable and $\phi^\e$ be a solution to \eqref{Benardgene}. 
 Then there exists a positive  constant
$C:=C(\nu, \kappa,K,L, T,M,N,\e_0)$ such that for any  $h\in {\mathcal A}_M$, $\e\in [0, \e_0]$,
\begin{equation}  \label{time}
I_n(h,\e):=\EX\Big[ 1_{G_N(T)}\; \int_0^T  |\phi_h^\e(s)-\phi_h^\e(\bar{s}_n)|^2 \, ds\Big]
\leq C\, 2^{-\frac{n}{2}}.
\end{equation}
\end{lemma}
\begin{proof}
Let $h\in {\mathcal A}_M$, $\e\geq 0$;  It\^o's formula yields $I_n(h,\e)=\sum_{1\leq i\leq 6} I_{n,i}$, where
\begin{align*}
I_{n,1}=&2\sqrt{\e}\;  \EX\Big( 1_{G_N(T)} \int_0^T \!\! ds \int_s^{\bar{s}_n}\!  \big(
\sigma(\phi_h^\e(r)) dW_r\, , \, \phi_h^\e(r)-\phi_h^\e(s) \big)\Big) , \\
I_{n,2}=&{\e}\;  \EX \Big( 1_{G_N(T)} \int_0^T  \!\!ds \int_s^{\bar{s}_n} \!\!
|\sigma(\phi_h^\e(r))|_{L_Q}^2 \, dr\Big) , \\
I_{n,3}=&2  \EX \Big( 1_{G_N(T)} \int_0^T \!\! ds \int_s^{\bar{s}_n} \!\! \big(
\tilde{\sigma}(\phi_h^\e(r)) \, h(r)\,  , \, \phi_h^\e(r)-\phi_h^\e(s) \big)\, dr\Big), \\
I_{n,4}=&- 2  \EX \Big( 1_{G_N(T)} \int_0^T  \!\! ds \int_s^{\bar{s}_n} \!\!
\big[ \nu \big(A_1 u_h^\e(r) , u_h^\e(r)-u_h^\e(s)\big)
+ \kappa \big(A_2 \theta_h^\e(r) ,  \theta_h^\e(r)-\theta_h^\e(s)\big)\big]dr\Big) ,  \\
I_{n,5}=&-2\,  \EX \Big( 1_{G_N(T)} \int_0^T  \!\! ds \int_s^{\bar{s}_n} \!\!
  \big(B( \phi_h^\e(r)),   \phi_h^\e(r)-\phi_h^\e(s)\big)
dr\Big) ,  \\
I_{n,6}=&2  \EX \Big( 1_{G_N(T)} \int_0^T \!\!  ds \int_s^{\bar{s}_n} \!\!
\big[  \big(u_{h,2}^\e(r),  \theta_h^\e(r)-\theta_h^\e(s)\big)
+  \big( \theta_h^\e(r) ,  u_{h,2}^\e(r)-u_{h,2}^\e(s)\big)\big]dr\Big) .
\end{align*}
Clearly  $G_N(T)\subset G_N(r)$ for $r\in [0,T]$. The
Burkholder-Davis-Gundy inequality, ({\bf A.4}) and the definition of
$G_N(r)$ yield for  $0\leq \e \leq \e_0$
\begin{eqnarray} \label{In1}
|I_{n,1}|&\leq & 2\sqrt{\e} \; \EX\int_0^T ds \Big| \int_s^{\bar{s}_n}
\big(
\sigma(\phi_h^\e(r)) dW_r\, , \, \phi_h^\e(r)-\phi_h^\e(s) \big)1_{G_N(r)} \Big| \nonumber\\
&\leq & 6\sqrt{\e} \int_0^T ds \; \EX \Big( \int_s^{\bar{s}_n} |\s(\phi_h^\e(r))|_{L_Q}^2
| \phi_h^\e(r)-\phi_h^\e(s)|^2 1_{G_N(r)} dr \Big)^{\frac{1}{2}} \nonumber \\
&\leq & 12 \sqrt{\e} \, \sqrt{KN(1+N)} \int_0^T ds \;  (T2^{-n})^{\frac{1}{2}} \leq C(\e_0,K,N,T) 2^{-\frac{n}{2}}.
\end{eqnarray}
Property ({\bf A.4}) implies that for $\e\leq \e_0$,
\begin{equation} \label{In2}
|I_{n,2}|\leq \e K \EX \Big( 1_{G_N(T)}   \int_0^T \!\! ds \int_s^{\bar{s}_n } \!\!
\big(1+|\phi_h^\e(r)|^2\big) dr\Big) \leq \e_0  K (1+N) T^2\, 2^{-n}.
\end{equation}
Schwarz's inequality, Fubini's theorem,   ({\bf A.4}) and the definition of ${\mathcal A}_M$ yield
\begin{eqnarray}  \label{In3}
|I_{n,3}|&\leq& 2 \sqrt{K}  \;  \EX \Big( 1_{G_N(T)} \int_0^T \!\! ds \int_s^{\bar{s}_n} \!\!
 \big(1+|\phi_h^\e(r)|^2\big)^{\frac{1}{2}}\, |h(r)|_0 |\, \phi_h^\e(r)-\phi_h^\e(s)|\, dr\Big)
\nonumber \\
& \leq&  4 \sqrt{KN(1+N)}\;  \EX \int_0^T |h(r)|_0\, dr \int_{\underline{r}_n}^r ds \leq C(K,N,M,T) 2^{-n}.
\end{eqnarray}
Schwarz's inequality and \eqref{boundgeneral}  imply that for some constant $\tilde{C}:=C(\e_0, \nu,\kappa,K,T)$
\begin{align} \label{In4}
I_{n,4} &\leq    \EX \Big( 1_{G_N(T)} \int_0^T  \!\! \! ds \int_s^{\bar{s}_n} \!\!\!
dr  \big[ - \nu \|u_h^\e(r)\|^2 - \kappa \|\theta_h^\e(r)\|^2 +
 \nu \|u_h^\e(r)\| \|u_h^\e(s)\| \nonumber \\
&\qquad  + \kappa \|\theta_h^\e(r)\| \|\theta_h^\e(s)\|\big]\Big) \nonumber \\
&\leq
 \frac{\nu+\kappa}{2}\;   \EX \Big( 1_{G_N(T)} \int_0^T ds \|\phi^\e_h(s)\|^2 \,
\int_s^{\bar{s}_n} dr \Big) \leq    \tilde{C}  \;   2^{-n}.
\end{align}
Inequalities \eqref{boundgeneral},  \eqref{normB1V'} and \eqref{inegB2},  Schwarz's inequality
and Fubini's theorem  imply that for some constant $\tilde{C}:=C(\e_0, \nu,\kappa,K,T)$, 
\begin{align} \label{In5}
|I_{n,5}| &\leq  2c_1\,  \EX \Big( 1_{G_N(T)} \int_0^T  \!\! ds \int_s^{\bar{s}_n} \!\!\!
dr \Big[ |u^\e_h(r)| \|u^\e_h(r)\|  \big( \|u_h^\e(r)\| + \| u^\e_h(s)\|\big]\big) \nonumber \\
&\qquad  +
|\phi_h^\e(r)|  \|\phi_h^\e(r)\| \big( \|\theta_h^\e(r)\| + \|\theta_h^\e(s)\|\big) \Big]
\nonumber \\
&
\leq 3\, c_1\, \sqrt{N} \EX\int_0^T dr \big( \|u^\e_h(r)\|^2
+ \|\phi^\e_h(r)\|^2\big) \int_{\underline{r}_n}^r ds
\nonumber \\
&\qquad + c_1\, \sqrt{N} \EX \int_0^T ds   \big( \|u^\e_h(s)\|^2
+ \|\phi^\e_h(s)\|^2\big) \int_s^{\bar{s}_n} dr
\leq  \sqrt{N} \tilde{C} 2^{-n} .
\end{align}
Finally, Schwarz's inequality implies that
\begin{equation} \label{In6}
|I_{n,6}|\leq 4 \EX \Big( 1_{G_N(T)} \int_0^T \!\! ds \int_s^{\bar{s}_n} \!\!\!
\big( |u_h^\e(r)|+ |u_h^\e(s)|\big)  \big( |\theta_h^\e(r)|+|\theta_h^\e(s)|\big) dr\Big)
\leq \frac{16 T^2 N}{ 2^n}.
\end{equation}
Collecting the upper estimates from \eqref{In1}-\eqref{In6}, we conclude the proof of \eqref{time}.
\end{proof}

Let $\e_0$ be defined as in Theorem \ref{wellposeness}
and  $(h_\e , 0 < \e \leq \e_0)$  be a family of random elements taking values in ${\mathcal A}_M$.
Let 
 $\phi^\e_{h_\e}$ be
the solution of the corresponding stochastic control equation
 with initial condition $\phi^\e_{h_\e}(0)=\xi \in H$:
\begin{eqnarray} \label{scontrol}
d \phi^\e_{h_\e} + [A\phi^\e_{h_\e} +B(\phi^\e_{h_\e})+R\phi^\e_{h_\e}]dt
=\s(\phi^\e_{h_\e}) h_\e dt+\sqrt{\e} \; \s(\phi^\e_{h_\e}) dW(t) . 
\end{eqnarray}
Note that $\phi^\e_{h_\e}={\mathcal G}^\e\Big(\sqrt{\e} \big( W_. + \frac{1}{\sqrt \e}
  \int_0^. h_\e(s)ds\big) \Big)$
due to the  uniqueness of the solution. \\

For all $\omega$ and  $h \in L^2([0, T], H_0)$, let $\phi_h$ be the solution of
the corresponding control equation (\ref{dc}) with initial condition
$\phi_h(0)=\xi(\omega)$:
\begin{eqnarray} \label{dcontrol2}
d \phi_h + [A\phi_h +B(\phi_h)+R\phi_h]dt =\s(\phi_h) h dt . 
\end{eqnarray}
Note that here we may assume that $h$ and $\xi$ are random,
but $\phi_h$ may defined pointwise by \eqref{dc}.

Let ${\mathcal C}_0=\{ \int_0^. h(s)ds \, :\, h\in L^2([0,T], H_0)\}  \subset C([0, T], H_0)$.
For every $\omega\in \Omega$, define ${\mathcal G}^0:   C([0, T], H_0)   \to X$ by
$ {\mathcal G}^0(g)(\omega)=\phi_h(\omega) $ for $  g=\int_0^. h(s)ds \in {\mathcal C}_0$
and ${\mathcal G}^0(g)=0$ otherwise.
\begin{prop} (Weak convergence) \label{weakconv}\\
 Suppose that $\sigma$ does not depend on time and satisfies Assumptions ({\bf A.1}), ({\bf A.4}) and ({\bf A.5}).
Let $\xi \in H$,   be ${\mathcal F}_0$-measurable such that $E|\xi|_H^4<+\infty$,
 and let $h_\e$ converge to $h$ in distribution as random elements
taking values in ${\mathcal A}_M$.  (Note that here ${\mathcal A}_M$
 is endowed with the weak topology
 induced by  the norm  defined in   \eqref{norm}).
Then
 as $\e \to 0$, $\phi^\e_{h_\e}$ converges in
distribution to $\phi_h$ in $X=C([0, T]; H) \cap L^2((0, T); V)$
endowed with  norm (\ref{norm}).
 That is,
  ${\mathcal G}^\e\Big(\sqrt{\e} \big( W_. + \frac{1}{\sqrt{\e}} \int_0^. h_\e(s)ds\big) \Big)$  converges in
distribution to $ {\mathcal G}^0\big(\int_0^. h(s)ds\big)$ in $X$, as $\e \to 0$.
\end{prop}
\begin{proof}
Since ${\mathcal A}_M$ is a Polish space (complete separable metric space),
by the Skorokhod representation theorem, we can construct
processes $(\tilde{h}_\e, \tilde{h}, \tilde{W})$ such that the
joint distribution of $(\tilde{h}_\e,  \tilde{W})$ is the same as
that of $(h_\e, W)$,  the distribution of $\tilde{h}$ coincides
with that of $h$, and $ \tilde{h}_\e  \to \tilde{h}$,  a.s., in
the (weak) topology of $S_M$. Hence a.s. for every $t\in [0,T]$, $\int_0^t \tilde{h}_\e(s) ds - \int_0^t
\tilde{h}(s)ds \to 0$ weakly in $H_0$.
Let $\Phi_\e=\phi^\e_{h_\e}-\phi_h$, or in component form
$\Phi_\e=(U_\e, \Theta_\e)=(u^\e_{h_\e}-u_h, \th^\e_{h_\e}-\th_h)$; then
\begin{eqnarray}\label{difference2}
& &d \Phi_\e + \big[A\Phi_\e +B(\phi^\e_{h_\e})-B(\phi_h)+R\Phi_\e\big]dt
\nonumber\\
&&\qquad  =\big[\s(\phi^\e_{h_\e}) h_\e -\s(\phi_h) h\big] dt +\sqrt{\e} \;
\s(\phi^\e_{h_\e}) dW(t), \;\; \Phi_\e(0)=0.
\end{eqnarray}
Let $\e_0$ be defined as in Theorem \ref{wellposeness}.
Set $\s_1=\s$, $\s_2=0$, $\bar{\s}=\s$, $\bar{\s}_1=0$, $\bar{\s}_2(s)=\s(\phi_h(s))\big(h_\e(s)-h(s)\big)$
and $\rho=0$.
Then $\phi_1=\phi^\e_{h_\e}$ and $\phi_2=\phi_h$ satisfy \eqref{phii}. Thus, Lemma \ref{Ito}, ({\bf A.4}) and ({\bf A.5})
yield for $0\leq \e\leq \e_0\wedge \frac{\nu \wedge \kappa}{4L}$:
\begin{align}  \label{total-error}
  |\Phi_\e(t)|^2 +&  (\nu\wedge \kappa) \int_0^t  \|\Phi_\e(s)\|^2 \, ds \leq   \sum_{i=1}^3 T_i(t,\e)
\nonumber  \\
\quad & + \int_0^t |\Phi_\e(s)|^2\, \Big[ 2+\frac{8 c_1^2}{\nu\wedge \kappa}\|\phi_h(s)\|^2 +
\frac{2L c_1 c_2}{\nu\wedge\kappa}\, |h(s)|_0^2\Big]\, ds,
\end{align}
where
\begin{align*}
T_1(t,\e)=& 2\sqrt{\e}\int_0^t \big( \Phi_\e(s), \s(\phi^\e_{h_\e}(s))\,  dW(s)  \big)\\
T_2(t,\e)= & \e K  \int_0^t (1+|\phi^\e_{h_\e}(s)|^2) ds , \\
T_3(t,\e)=& 2\int_0^t \big( \s(\phi_h(s))\, \big( h_{\e}(s)-h(s)\big) ,  \Phi_\e(s)\big)\, ds.
\end{align*}
Our goal here is to show that as $\e \to 0$,  $\sup_{0\leq t\leq
T} |\Phi_\e(t)|^2 + \int_0^T \|\Phi_\e(s)\|^2ds \to 0$ in
probability, which implies
 that $\phi_{h_\e} \to \phi_h$  in distribution in
$X:=C([0, T]; H) \cap L^2((0, T); V)$.
 Fix $N>0$ and for $t\in [0,T]$ let
\begin{eqnarray*}  G_N(t)&=&\Big\{ \sup_{0\leq s\leq t} |\phi_h(s)|^2 \leq N\Big\} \cap
 \Big\{ \int_0^t \|\phi_h(s)\|^2  ds \leq N
\Big\} , \\
G_{N,\e}(t)&=&  G_N(t)\cap \Big\{ \sup_{0\leq s\leq t} |\phi^\e_{h_\e}(s)|^2 \leq N\Big\} \cap
 \Big\{ \int_0^t \|\phi^\e_{h_\e}(s)\|^2  ds \leq N
\Big\} .
\end{eqnarray*}
\emph{\textbf{Claim 1.}}
For any $\e_0>0$, $ {\displaystyle \sup_{0<\e\leq \e_0}\; \sup_{h,h_\e \in {\mathcal A}_M}
\PX(G_{N,\e}(T)^c )\to 0 \;
\mbox{\rm as }\; N\to \infty.}$
\\
Indeed, for $\e>0$, $h,h_\e \in {\mathcal A}_M$, the Markov inequality and
 estimate \eqref{boundgeneral} imply
\begin{align*}
&  \PX ( G_{N,\e}(T)^c )
\leq   \PX \Big(\sup_{0\leq s\leq T} |\phi_h(s)|^2 > N \Big)
 +\PX \Big(\sup_{0\leq
s\leq T} |\phi^\e_{h_\e}(s)|^2 > N \Big) \\
&\qquad \qquad\qquad   +  \PX\  \Big( \int_0^T \big(\|\phi_h(s)\|^2 ds  >N\Big)
+ \PX\Big(\int_0^T \|\phi^\e_{h_\e}(s)\|^2 \big)ds > N \Big)\\
&\leq  \frac{1}{N}
  \sup_{ 
\;h, h_\e \in {\mathcal A}_M}
 \EX
\Big( \sup_{0\leq s\leq T} |\phi_h(s)|^2 + \sup_{0\leq s\leq T}
|\phi^\e_{h_\e}(s)|^2  
+ \int_0^T(\|\phi_h(s)\|^2+ \|\phi^\e_{h_\e}(s)\|^2 )ds \Big) \\
&\leq
 {C_1(\nu, \k,K,L, T, M) \, \big(1+ E|\xi|^4\big)}{N}^{-1}.
\end{align*}
\emph{\textbf{Claim 2.}} For fixed $N>0$,
 $h, h_\e \in {\mathcal A}_M$  such that  as $\e \to 0$, $h_\e\to h$ a.s. in the weak topology
 of $L^2([0,T],H_0)$,
 one has  as $\e \to 0$
\begin{equation} \label{cv1}
\EX\Big[ 1_{G_{N,\e}(T)} \Big( \sup_{0\leq t\leq T } |\Phi_\e(t)|^2 + \int_0^T  \|\Phi_\e(t)\|^2 \, dt\Big)
\Big] \to 0.
\end{equation}
Indeed,  \eqref{total-error}  and Gronwall's lemma imply that on $G_{N,\e}(T)$,
\[ \sup_{0\leq t\leq T} |\Phi_\e(t)|^2 \leq \Big[ \sup_{0\leq t\leq T}
 \big( T_1(t,\e) + T_3(t,\e)\big) +
\e KT(1+N)\Big] \, 
e^{ 2T+\frac{8 \,c_1^2\, N}{\nu\wedge\kappa}+\frac{2\,L\,c_1\, c_2\,M}{\nu\wedge \kappa}  }
.\]
Thus, using again \eqref{total-error}  we deduce that for some constant $\tilde{C}=C(\nu,\kappa,K,L, T,M,N)$,
 one has for every $\e>0$:
\begin{equation} \label{*}
\EX\big( 1_{G_{N,\e}(T)} \, |\Phi_\e|_X^2 \big) \leq \tilde{C} \Big( \e KT(1+N) + \EX \Big[ 1_{G_{N,\e}(T)} \,
\sup_{0\leq t\leq T} \big( T_1(t,\e) + T_3(t,\e)\big) \Big] \Big).
\end{equation}
Since the sets $G_{N,\e}(.)$ decrease, $\EX\big(1_{G_{N,\e}(T)}
 \sup_{0\leq t\leq T} |T_1(t,\e)|\big) \leq
\EX(\lambda_\e)$, where
\[ \lambda_\e := 2\sqrt{\e} \sup_{0\leq t\leq T }
\Big|\int_0^t 1_{G_{N,\e}(s)} \big( \Phi_\e(s), \s(\phi^\e_{h_\e}(s)) dW(s) \big)\Big|.\]
The scalar-valued random variables $\lambda_\e$ converge to 0 in $L^1$ as $\e\to 0$.
Indeed, by the Burkholder-Davis-Gundy inequality, ({\bf A.4}) and the definition of $G_{N,\e}(s)$,
we have
\begin{eqnarray}
\EX (\lambda_\e) & \leq&  6\sqrt{\e} \; \EX
\Big\{\int_0^T 1_{G_{N,\e}(s)} \,  |\Phi_\e(s)|^2 \;
 |\s(\phi^\e_{h_\e}(s))|^2_{L_{Q}}
 ds\Big\}^\frac12 \nonumber \\
&\leq & 6\sqrt{\e} \; \EX \Big[  \Big\{ 4N \int_0^T 1_{G_{N,\e}(s)}\,
K\, (1+ | \phi^\e_{h_\e}(s) |^2 )  ds\Big\}^\frac12 \Big] \nonumber \\
& \leq & 12  \,  \sqrt{\e} \, \sqrt{KT} \, (1+N)   . \label{lambda1}
\end{eqnarray}
For $k=0, \cdots, 2^n$ set $t_k=kT2^{-n}$; for  $s\in ]t_k, t_{k+1}]$,  set $\bar{s}_n=t_{k+1}$ and $\underline{s}_n
= t_k$.  Then for any $n\geq 1$,
\[ \EX\Big( 1_{G_{N,\e}(T)}\sup_{0\leq t\leq T} |T_3(t,\e)| \Big)
\leq 2\;  \sum_{i=1}^3  \tilde{T}_i(N,n, \e)+ 2 \; \EX \big( \bar{T}_4(N,n,\e,\omega)\big),\]
 where
\begin{align*}
\tilde{T}_1(N,n,\e)=& \EX \Big[ 1_{G_{N,\e}(T)} \sup_{0\leq t\leq T}  \Big| \int_0^t  \Big(  \s(\phi_h(s))
  \big(h_\e(s)-h(s)\big) \, ,\,
\big[\Phi_\e(s)-\Phi_\e(\bar{s}_n)\big] \Big)  ds\Big| \Big] ,\\
\tilde{T}_2 (N,n,\e)=& \EX\Big[  1_{G_{N,\e}(T)} \sup_{0\leq t\leq T} \Big| \int_0^t
\Big( \big[ \s(\phi_h(s)) - \s(\phi_h(\bar{s}_n))\big]
\big(h_\e(s) - h(s) \big)\, ,\, \Phi_\e(\bar{s}_n)\Big) ds\Big| \Big] ,\\
\tilde{T}_3(N,n,\e)=& \EX \Big[  1_{G_{N,\e}(T)} \sup_{1\leq k\leq 2^n}  \sup_{t_{k-1}\leq t\leq t_k}
\Big| \big( \sigma(\phi_h(t_k))  \int_{t_{k-1}}^t (h_\e(s)-h(s))\, ds \, ,\, \Phi_\e(t_k)\big) \Big| \Big]\\
\bar{T}_4(N,n, \e)=&  1_{G_{N,\e}(T)} \sum_{k=1}^{2^n} \Big| \Big( \s(\phi_h(t_k))
\int_{t_{k-1}}^{t_k }  \big(h_\e(s)-h(s)\big)\, ds \, ,\,
\Phi_\e(t_k )\Big) \Big| .
\end{align*}
Using Schwarz's inequality, ({\bf A.4}) and Lemma \ref{timeincrement}, we deduce that  for some constant
$\bar{C}_1:= C(\nu,\kappa,K,T,M,N)$ and any   $\e \in ]0, \e_0]$,
\begin{align} \label{eqT1}
&  \tilde{T}_1(N,n,\e)\leq  \sqrt{K} \EX\Big[ 1_{G_{N,\e}(T)}  \int_0^T
 \big( 1+|\phi_h(s)|^2\big)^{\frac{1}{2}}
|h_\e(s)-h(s)|_0\, \big| \Phi_\e(s)-\Phi_\e(\bar{s}_n)\big|\, ds\Big]  \nonumber \\
& \quad   \leq  \sqrt{2 K (1+N)}\;  \Big( \EX  \int_0^T |h_\e(s)-h(s)|_0^2\, ds \Big)^{\frac{1}{2}}
\nonumber \\
&\qquad \times \Big(  \EX \Big[ 1_{G_{N,\e}(T)}  \int_0^T
 \big\{ |\phi^\e_{h_\e}(s) - \phi^\e_{h_\e}(\bar{s}_n)|^2 +
 |\phi_{h}(s) - \phi_{h}(\bar{s}_n)|^2 \big\}\, ds\Big] \Big)^{\frac{1}{2}}
\nonumber \\
 &\quad
\leq \bar{C}_1 \;  2^{-\frac{n}{4}}.
\end{align}
A similar computation based on ({\bf A.5}) and Lemma \ref{timeincrement} yields  for
some constant $\bar{C}_2:=C(\nu,\kappa,K,L, T,M,N)$ and any   $\e \in ]0, \e_0]$
\begin{align} \label{eqT2}
 \tilde{ T}_2(N,n,\e) &\leq  \sqrt{L} \Big( \EX \Big[ 1_{G_{N,\e}(T)}  \int_0^T\!\!
 |\phi_{h}(s) - \phi_{h}(\bar{s}_n)|^2 \, ds\Big]  \Big)^{\frac{1}{2}} \Big(
 \EX \int_0^T \!\! |h_\e(s)-h(s)|_0^2 \, 4 N  \, ds \Big)^{\frac{1}{2}}
 \nonumber \\
& \leq \bar{C}_2 \;  2^{-\frac{n}{4}}.
\end{align}
Using Schwarz's inequality and ({\bf A.4}) we deduce  for $\bar{C}_3=C(K,N,M)$ and any $\e \in ]0, \e_0]$
\begin{align} \label{eqT3}
\tilde{T}_3(N,n,\e)&\leq \sqrt{K} \EX \Big[  1_{G_{N,\e}(T)} \sup_{1\leq k\leq 2^n} \big(1+| \phi_h(t_k)|^2
\big)^{\frac{1}{2}} \int_{t_{k-1}}^{t_k}\!\! |h_\e(s)-h(s)|_0 \, ds \, |\Phi_\e(t_k)| \Big]\nonumber \\
&\leq 2 \sqrt{KN(1+N)}\;  \EX\Big( \sup_{1\leq k\leq 2^n}  \int_{t_{k-1}}^{t_k} |h_\e(s)-h(s)|_0 \, ds\Big)
 \nonumber \\
&\leq 8  \sqrt{KN(1+N)} \; \sqrt{M} \;  2^{-\frac{n}{2}} = \bar{C}_3\;  2^{-\frac{n}{2}}.
\end{align}
Finally, note that the weak convergence of $h_\e$ to $h$ implies that for any $a,b\in [0,T]$, $a<b$,
as $\e\to 0$, the integral $\int_a^b h_\e(s) ds \to \int_a^b h(s) ds$ in the weak topology of $H_0$. Therefore, since for
$\phi\in H$ the operator $\sigma(\phi)$ is compact from $H_0$ to $H$, we deduce that
$\left| \sigma(\phi) \big(   \int_a^b h_\e(s) ds - \int_a^b h(s) ds  \big) \right|_H \to 0$ as $\e \to 0$.
Hence a.s. for fixed $n$ as $\e \to 0$, $\bar{T}_4(N,n,\e,\omega) \to 0$. Furthermore, $\bar{T}_4(N,n,\e,\omega)
\leq \sqrt{K} \sqrt{1+N} \sqrt{4N} \sqrt{M}$ and hence the dominated convergence theorem proves that for any
fixed $n$, $\EX(\bar{T}_4(N,n,\e))\to 0$ as $\e \to 0$.

Thus, given $\alpha >0$, we may choose $n_0$ large enough to have
 $(\bar{C}_1 + \bar{C}_2) 2^{-\frac{n}{4}} + \bar{C}_3 2^{-\frac{n}{2}}\leq \alpha$ for $n\geq n_0$.
 Then for fixed
$n\geq n_0$, let  $\e_1\in ]0,\e_0]$  be  such that for $0<\e\leq \e_1$, $\EX \big[\bar{T}_4(N,n,\e)\big]
\leq \alpha$. Using \eqref{eqT1}-\eqref{eqT3}, we deduce that for $\e\in ]0,\e_1]$,
\begin{equation} \label{eqT4}
\EX \Big[ 1_{G_{N,\e}(T)} \sup_{0\leq t\leq T} |T_3(t,\e)|\Big] \leq 2\alpha.
\end{equation}
Claim 2 is a straightforward consequence of  inequalities \eqref{*}, \eqref{lambda1} and \eqref{eqT4}.
\smallskip

To conclude the proof of  Proposition \ref{weakconv}, let
$\delta>0$ and $\alpha >0$ and set
\[  \Lambda_\e := |\Phi_\e|^2_X =  \sup_{0\leq t\leq T}  |\Phi_\e(t)|^2
+ \int_0^T \|\Phi_\e(s)\|^2 ds.   \]
Then the Markov inequality implies that
\[
\PP(\Lambda_\e > \de )  =  \PP(G_{N,\e}(T)^c )+ \frac{1}{\delta} \EX\Big( 1_{G_{N,\e}(T)} |\Phi_\e|^2_X\Big)
\]
Using \emph{Claim 1}, one can choose $N$ large enough to
make sure that $ \PP(G_{N,\e}(T)^c)<\alpha$ for every $\e\leq \e_0$. Fix $N$~;
{\it Claim 2} shows that for  $\e$ small enough,
 $ \EX\Big( 1_{G_{N,\e}(T)} |\Phi_\e|^2_X\Big)
< \delta  \alpha$. This concludes the proof of the proposition.
\end{proof}

The following compactness result
will show that the rate function of the LDP  satisfied by the solution
 to \eqref{scontrol} is a good rate function. The proof is similar to that of Proposition \ref{weakconv}
and easier.
\begin{prop}  \label{compact} (Compactness)\\
Let $M$ be any fixed finite positive number and let $\xi\in H$ be deterministic. Define
$$
K_M=\{\phi_h \in  C([0, T]; H) \cap L^2((0, T); V):  h \in S_M \},
$$
where $\phi_h$ is the unique solution of the deterministic control
equation:
\begin{eqnarray} \label{dcontrol}
d \phi_h(t) + \big[A\phi_h(t) +B(\phi_h(t))+R\phi_h(t)\big]dt =\s(\phi_h(t)) h(t) dt, \;\;
\phi_h(0)=\xi,
\end{eqnarray}
and $\s$ does not depend on time and satisfies ({\bf A.1}), ({\bf A.4}) and ({\bf A.5}).
Then $K_M$ is a compact subset of  $ X$.
\end{prop}
\begin{proof}
Let $(\phi_n)$ be a sequence in $K_M$, corresponding to solutions of
(\ref{dcontrol}) with controls $(h_n)$ in $S_M$:
\begin{eqnarray} \label{dcontroln}
d \phi_n(t) + \big[A\phi_n(t) +B(\phi_n(t))+R\phi_n(t)\big]dt =\s(\phi_n(t)) h_n(t) dt, \;\;
\phi_n(0)=\xi.
\end{eqnarray}
 Since $S_M$ is a
bounded closed subset in the Hilbert space $L^2((0, T); H_0)$, it
is weakly compact. So there exists a subsequence of $(h_n)$, still
denoted as $(h_n)$, which converges weakly to a limit $h$ in
$L^2((0, T); H_0)$. Note that in fact $h \in S_M$ as $S_M$ is
closed. We now show that the corresponding subsequence of
solutions, still denoted as  $(\phi_n)$, converges in $ X$
 to $\phi$ which is the solution of the
following ``limit'' equation
\begin{eqnarray} \label{limiteqn}
d \phi(t) + [A\phi(t) +B(\phi(t))+R\phi(t)]dt =\s(\phi(t)) h(t) dt, \;\;
\phi(0)=\xi.
\end{eqnarray}
This will complete the proof of the compactness of $K_M$.
  To ease notation
 we will often drop the time parameters $s$, $t$, ... in the equations and integrals.

Let $\Phi_n=\phi_n-\phi$, or in component form $\Phi_n=(U_n,
\Theta_n)=(u_n-u, \th_n-\th)$; then
\begin{eqnarray}\label{difference}
d \Phi_n + [A\Phi_n +B(\phi_n)-B(\phi)+R\Phi_n]dt =[\s(\phi_n) h_n
-\s(\phi) h] dt, \;\; \Phi_n(0)=0.
\end{eqnarray}
Set $\s_1=\s_2=0$, $\bar{\s}=\s$, $\bar{\s}_1=0$, $\bar{\s}_2(s)=\s(\phi(s))\, [h(s)-h_n(s)]$,
$h_\e=h_n$, $\rho=0$. Then $\phi_1:=\phi_n$ and $\phi_2:=\phi$ satisfy \eqref{phii}.

Thus Lemma \ref{Ito} yields the following  integral
inequality
\begin{align} \label{error}
 |\Phi_n(t)|^2 + &(\nu\wedge \k)\, \int_0^t \|\Phi_n(s)\|^2ds
  \leq 2 \int_0^t \big( \s(\phi(s))\, [h(s)-h_n(s)]\, ,\, \Phi_n(s)\big) \, ds  \nonumber \\
 & +  \int_0^t
 \Big\{2+ \frac{8 c_1^2}{\nu\wedge\kappa} \|\phi(s)\|^2 + \frac{2Lc_1c_2}{\nu\wedge\kappa}\,   |h_n(s)|_0^2
\Big\}  \,  |\Phi_n(s)|^2 ds.
\end{align}
For $N\geq 1 $  and $k=0,\cdots,2^N$, set $t_k=k2^{-N}$. For $s\in ]t_{k-1}, t_k]$, $1\leq k\leq 2^N$,
let $\bar{s}_N=t_k$.
Inequality \eqref{control-norm2} implies that there exists a constant $\bar{C}>0$ such that
 \[ \sup_n \Big[ \sup_{0\leq t\leq T} \big(|\phi(t)|^2 + |\phi_n(t)|^2\big) + \int_0^T
\big( \|\phi(s)\|^2+\|\phi_n(s)\|^2\big)ds \Big] =  \bar{C} <+\infty.\]
 Thus  Gronwall's inequality implies
\begin{equation} \label{errorbound}
 \sup_{ t\leq T} |\Phi_n(t)|^2  \leq
\exp\Big( 2T + \frac{8c_1^2\bar{C} }{\nu\wedge\kappa} +\frac{2Lc_1c_1 M}{\nu\wedge\kappa} \Big)
  \sum_{i=1}^4 I_{n,N}^{i}  ,
\end{equation}
where
\begin{eqnarray*}
I_{n,N}^1&=& \int_0^T \big| \big( \s(\phi(s))\,  [h_n(s)- h(s)]\, ,\, \Phi_n(s)-\Phi_n(\bar{s}_N)\big)\big|\, ds,  \\
I_{n,N}^2&=& \int_0^T\Big|  \Big( \big[ \s(\phi(s)) - \s(\phi(\bar{s}_N))\big]  [h_n(s)-h(s)]\, ,\,
\Phi_n(\bar{s}_N)\Big)\Big| \, ds,  \\
I_{n,N}^3&=& \sup_{1\leq k\leq 2^N} \sup_{t_{k-1}\leq t\leq t_k} \Big|\Big( \s(\phi(t_k ))
 \int_{t_{k-1 }}^t (h_n(s)-h(s))
ds \; ,\; \Phi_n(t_k) \Big)\Big| , \\
I_{n,N}^4&=& \Big| \sum_{k=1}^{2^N} \Big( \s(\phi(t_k )) \,  \int_{t_{k-1}}^{t_k }
[ h_n(s)-h(s)]\, ds   \; ,\; \Phi_n(t_k ) \Big) \Big| .
\end{eqnarray*}
Schwarz's inequality, ({\bf A.4}), ({\bf A.5}) and Lemma \ref{timeincrement}
 imply that for some constant $C$ which does not
depend on $n$ and $N$,
\begin{align} \label{estim1}
I_{n,N}^1 &\leq  \Big( \int_0^T\!\! K (1+\bar{C}) |h_n(s)-h(s)|_0^2 ds \Big)^{\frac{1}{2}} \nonumber\\
&\quad\times \Big( 2 \int_0^T \!\! \big( |\phi_n(s)-\phi_n(\bar{s}_N)|^2 + |\phi(s)-\phi(\bar{s}_N)|^2\big) ds
 \Big)^{\frac{1}{2}}
\nonumber\\
&\leq C 2^{-\frac{N}{4}} \, ,\\
I_{n,N}^2&\leq  \Big( L\int_0^T |\phi(s)-\phi(\bar{s}_N)|^2 ds \Big)^{\frac{1}{2}}
 \Big( \bar{C} \int_0^T |h_n(s)-h(s)|_0^2\, ds\Big)^{\frac{1}{2}} \leq C 2^{-\frac{N}{4}}\, ,
\label{estim2}\\
I_{n,N}^3&\leq K \big( 1+\sup_t |\phi(t)|) \sup_t\big( |\phi(t)|+\phi_n(t)|\big) 2^{-\frac{N}{2}} 2M \leq
C \, 2^{-\frac{N}{2}}\,  .
\label{estim3}
\end{align}
Thus, given $\alpha >0$, one may choose $N$ large enough to have $ \sup_n \sum_{i=1}^3 I_{n,N}^i
\leq \alpha$. Then, for fixed $N$ and $k=1, \cdots, 2^N$,
 as $n\to \infty$, the weak convergence of $h_n$ to $h$ implies
that of $\int_{t_{k-1}}^{t_k} (h_n(s)-h(s))ds$ to 0 weakly in $H_0$. Since $\s(\phi(t_k))$ is a compact operator,
we deduce that for fixed $k$ the sequence $\s(\phi(t_k)) \int_{t_{k-1}}^{t_k} (h_n(s)-h(s))ds$
converges to 0 strongly in $H$ as $n\to \infty$.
Since $\sup_n \sup_k |\Phi_n(t_k)|\leq 2\tilde{C}$, we have
 $\lim_n I_{n,N}^4=0$. Thus as $ n \to \infty$,
$
\sup_{0\leq t\leq T}|\Phi_n(t)|^2 \to 0. 
$
Using this convergence and \eqref{errorbound}, we deduce that $\|\Phi_n\|_X \to 0$ as $n\to \infty$.
This shows that every sequence in $K_M$ has a convergent
subsequence. Hence $K_M$ is  a compact subset of $X$.
\end{proof}
With the above results, we have the following large deviation
theorem.
\begin{theorem}  
 Suppose that  $\s$ does not depend on time and
satisfies  ({\bf A.1}), ({\bf A.4}) and ({\bf A.5}), 
 let  $\phi^\e$ be the solution of the stochastic B\'{e}nard problem
\eqref{Benard}.
Then $\{\phi^\e\}$ satisfies the large deviation principle in
$C([0, T]; H) \cap L^2((0, T); V)$,  with the good rate function
\begin{eqnarray} \label{ratefc}
 I_\xi (\psi)= \inf_{\{h \in L^2(0, T; H_0): \; \psi ={\mathcal G}^0(\int_0^. h(s)ds) \}}
 \Big\{\frac12 \int_0^T |h(s)|_0^2\,  ds \Big\}.
\end{eqnarray}
Here the infimum of an empty set is taken as infinity.
\end{theorem}
\begin{proof}
 Propositions \ref{compact} and  \ref{weakconv} imply that
$\{\phi^\e\}$ satisfies the Laplace principle which is equivalent
to the large deviation principle in $X=C([0, T], H) \cap L^2((0,
T), V)$ with the above-mentioned rate function;
 see Theorem 4.4 in \cite{BD00} or
Theorem 5 in \cite{BD07}.
\end{proof}
\medskip

\noindent {\bf Acknowledgements.} This work was partially done
while J. Duan was a  Professor Invit\'{e} at the Universit\'{e}
Paris 1, France. J. Duan would like to thank  this University for
financial support and the colleagues at Samos-Matisse, Centre
d'\'Economie de la Sorbonne, for their hospitality.   It was
completed while A. Millet was visiting the Mittag Leffler
Institute, Sweden,   which provided financial support. She would like to thank the center for
 excellent working conditions  and
 a very friendly  atmosphere.
The authors   would also like to thank Padma Sundar
  for   helpful discussions and an anonymous referee for his careful reading and helpful comments.


\end{document}